\documentclass[12pt]{article}

\usepackage{amsmath,amssymb,amsfonts,amsthm}

\newtheorem{theorem}{Theorem}[section]
\newtheorem{claim}[theorem]{Claim}
\newtheorem{lemma}[theorem]{Lemma}
\newtheorem{proposition}[theorem]{Proposition}
\newtheorem{corollary}[theorem]{Corollary}

\theoremstyle{definition}
\newtheorem{definition}[theorem]{Definition}
\newtheorem{example}[theorem]{Example}

\theoremstyle{remark}
\newtheorem{remark}[theorem]{Remark}

\newcommand{\ad}{\operatorname{ad}}

\newcommand{\Ann}{\operatorname{Ann}}

\newcommand{\Aut}{\operatorname{Aut}}

\newcommand{\Casim}{\operatorname{Casim}}

\newcommand{\const}{\operatorname{const}}
\newcommand{\Coup}{\operatorname{Coup}}

\newcommand{\Curv}{\operatorname{Curv}}

\newcommand{\End}{\operatorname{End}}

\renewcommand{\flat}{\operatorname{flat}}
\newcommand{\flow}{\operatorname{flow}}

\newcommand{\Ham}{\operatorname{Ham}}
\newcommand{\hor}{\operatorname{hor}}

\newcommand{\id}{\operatorname{id}}
\newcommand{\Inn}{\operatorname{Inn}}

\newcommand{\Poiss}{\operatorname{Poiss}}

\newcommand{\rank}{\operatorname{rank}}

\newcommand{\so}{\operatorname{so}}

\newcommand{\Sect}{\operatorname{Sect}}

\begin{document}

\title{Poisson Equivalence over a Symplectic Leaf}

\author{Yurii Vorobjev\thanks{Research was partially supported
by CONACYT, grant 43208, 
and by the Russian Foundation for Basic Research,
grant 02-01-00952.}\\
\small\it Department of Mathematics, University of Sonora,
Mexico; \\ 
\small\it Moscow Institute of Electronics and Mathematics,
Russia\\
\small yurimv@guaymas.uson.mx; vorob@miem.edu.ru}

\date{}

\maketitle

\begin{abstract}
We study the equivalence of Poisson structures around a given
symplectic leaf of nonzero dimension. 
Some criteria of Poisson equivalence are derived 
from a homotopy argument for coupling Poisson structures. 
In the case when the transverse Lie algebra of the symplectic
leaf is semisimple of compact type, we show that an obstruction
to the linearizability is the cohomology class of a Casimir
$2$-cocycle. This allows us to obtain a semilocal analog of
the Conn linearization theorem and to clarify examples of
nonlinearizable Poisson structures due to~\cite{DW}. 
\end{abstract}

\section{Introduction}

This work is devoted to the problem of classification of
Poisson structures near a given symplectic leaf of nonzero
dimension.

Let $(M,\Psi)$ be a Poisson manifold with Poisson tensor $\Psi$.
According to the local splitting theorem~\cite{We1}, in a
neighborhood of each point $m\in M$, 
$\Psi$ is equivalent to the direct product of a nondegenerate
Poisson structure and a transverse Poisson structure vanishing
at~$m$.  
The linearization of the transverse factor at $m$ leads to the
\textit{linearized transverse Poisson structure} $\Lambda_{m}$ 
which leaves naturally on the normal space $E_{m}$ to the
symplectic leaf through $m$. This structure is 
uniquely determined by the \textit{transverse Lie algebra}
$\mathfrak{g}_{m}$ of $m$. If $m$ is a singular point (a
zero-dimensional leaf), $\Psi(m)=0$, then 
$E_{m}=T_{m}M$ and the linearized transverse Poisson structure
$\Lambda_{m}$ gives a linear approximation to $\Psi$ at $m$. 
In this case, one can state the linearization problem~\cite{We1}
which consists of determining whether the original Poisson
structure is locally equivalent to its linear approximation. 
In formal, $C^{\infty}$, and analytic settings, this problem was
studied in \cite{Ar,We1,Co1,Co2,Du1}. In particular, 
it was shown in~\cite{We1} that if
the transverse Lie algebra $\mathfrak{g}_{m}$ is semisimple, then
there exists a formal linearization. Later, for analytic Poisson
structures, this fact was proved in~\cite{Co1}. 
In the smooth case, the local linearization theorem
due to~\cite{Co2} says that the Poisson structure $\Psi$
vanishing at $m$ is locally isomorphic to the linear
approximation $\Lambda_{m}$ if $\mathfrak{g}_{m}$ 
is a semisimple Lie algebra of compact type. Further
developments and generalizations of these results where obtained
in \cite{We3,Du2,FeM,DZ,Zu1,Zu2}. A review of recent results on
the local linearization problem can be found in~\cite{FeM}.

We are interested in the linearization problem, 
or more generally, normal forms for the Poisson structure $\Psi$
in the semilocal context, around a symplectic leaf~$B$. 
A more interesting situation occurs when 
the symplectic foliation nearby~$B$
is singular. The Poisson topology of neighborhoods of (singular)
symplectic leaves has been studied in several papers
\cite{GiGo,Fe1,Fe2}.  
The recent work~\cite{CrFe} devoted is to rigidity and flexibility
phenomena in Poisson geometry.

Our approach is based on the notion of a \textit{coupling
Poisson structure\/} \cite{Vo1,Va2}. 
The coupling procedure for symplectic structures was
introduced in~\cite{St} and further developed in~\cite{GLS}. The
Poisson coupling can be defined on fibered spaces~\cite{Vo1} and
foliated manifolds~\cite{Va2} as well. 
The key observation is that a given Poisson structure
$\Psi$ can be realized as a coupling Poisson structure on the
normal bundle $E=TM/TB$ of the leaf. 
Using a diffeomorphism $\mathbf{f}:\,E\rightarrow M$ identical
on~$B$, one can move from $M$ to the total space~$E$.
The result is the
Poisson tensor $\Pi=\mathbf{f}^{\ast}\Psi$ on $E$ with the same
symplectic leaf $B$ (as the zero section). 
Then, 
in a neighborhood of $B$ in $E$ the Poisson tensor $\Pi$ induces
an \textit{intrinsic Ehresmann connection} $\Gamma$
corresponding to the splitting
\begin{equation}
TE=\mathbb{H}\oplus\mathbb{V},
\tag{1.1}
\end{equation}
where $\mathbb{V}$ is the vertical subbundle of
$\pi:\,E\rightarrow B$ and the 
horizontal subbundle $\mathbb{H}$ is uniquely determined by $\Pi$,
\begin{equation} 
\mathbb{H}=\Pi^{\natural}(\Ann\mathbb{V}).
\tag{1.2}
\end{equation}
Here $\Pi^{\natural}:\,T^{\ast}E\rightarrow TE$ 
is the bundle morphism induced by $\Pi$
$\Ann\mathbb{V}$ is the annihilator of $\mathbb{V}$.
Furthermore, the bivector field $\Pi$ has the following
decomposition with respect to~(1.1):
\begin{equation}
\Pi=\Pi_{H}+\Pi_{V},
\tag{1.3}
\end{equation}
where $\Pi_{H}$ and $\Pi_{V}$ are the horizontal and vertical
bivector fields, respectively. 
In general, the decomposition of $\Pi$ involves a bivector field
of degree $(1,1)$. The vertical part $\Pi_{V}$ is a Poisson
tensor with the property: 
the symplectic leaf of $\Pi_{V}$ through each point $m\in E$ belongs
to the fiber $E_{\xi}$ over $\xi=\pi(m)$. Such Poisson tensors on fibered
(or foliated) spaces are called fiber-tangent (or leaf-tangent) Poisson
structures~\cite{Va2}. 
The restriction of $\Pi_{V}$ to each fiber $E_{\xi}$
vanishes at $0$ and defines the transverse Poisson structure of
$\xi$. The horizontal part $\Pi_{H}$ is not a Poisson tensor in
general. The Jacobi identity for $\Pi_{H}$ is equivalent to the
zero curvature condition for $\Gamma$. 
Because of the property $\mathbb{H}|_{B}=TB$, in a neighborhood
of $B$, the horizontal bivector field $\Pi_{H}$ is uniquely
defined, by a nondegenerate $C^{\infty}(E)$-valued $2$-form $F$
on $B$, which is called a \textit{coupling form}. 
This form can be viewed as a ``perturbation'' of the
symplectic structure of the leaf $B$. Thus, the Poisson tensor
is defined by the triple $(\Pi_{V},\Gamma,F)$ and formula (1.3)
can be interpreted as a result of coupling the form $F$ and the
fiber-tangent Poison structure $\Pi_{V}$ via the connection
$\Gamma$. We call $\Pi$ a \textit{coupling Poisson tensor}
associated with geometric data $(\Pi_{V},\Gamma,F)$. The Jacobi
identity for $\Pi$ is equivalent to some ``integrability''
conditions for $(\Pi_{V},\Gamma,F)$ which have a natural
geometric sense. Thus,  given some integrable
geometric data we can reconstruct the coupling Poisson tensor.

To define a linearized Poisson structure of $\Pi$ at $B$, first,
one can linearize the geometric data
$(\Pi_{V},\Gamma,F)$ and then construct a ``new'' coupling
Poisson tensor. The main point here is that the linearization
procedure for $(\Pi_{V},\Gamma,F)$ preserves the integrability
conditions.  Thus, we get the linearized geometric data
$(\Lambda,\Gamma^{(1)},F^{(1)})$
consisting of a linearized transverse Poisson structure
$\Lambda$, a homogeneous (linear) connection $\Gamma^{(1)}$ 
on $E$ and a linearized coupling form $F^{(1)}$. 
We define
the \textit{linearized Poisson structure} of $\Pi$ at $B$ as the
coupling Poisson tensor $\Pi^{(1)}$ associated with the data
$(\Lambda,\Gamma^{(1)},F^{(1)})$. The Poisson tensor
$\Pi^{(1)}$ is independent of the choice of a diffeomorphism
$\mathbf{f}$ up to an isomorphism and hence defines an intrinsic
infinitesimal characteristic of the original Poisson structure
$\Psi$ at $B$ \cite{Vo1,Vo2}. Actually, as was shown in
\cite{Vo1,Vo3}, $\Pi^{(1)}$ is uniquely determined by the
\textit{transitive Lie algebroid\/} of the leaf $B$~\cite{We2}. 
A derivation of $\Pi^{(1)}$ in the context of the Lie
algebroid approach, can also be found in~\cite{Va2}.

Now, one can state the semilocal linearization problem saying
that the Poisson structure $\Psi$ is \textit{linearizable} at
$B$ if $\Pi$ and $\Pi^{(1)}$ are isomorphic by a diffeomorphism
$\phi:\,E\rightarrow E$ from an appropriate class.  We assume
that $\phi$ belongs to the group of diffeomorphisms on $E$
satisfying the conditions $\phi|_{B}=\id_{B}$  
and $d_{B}\phi| _{E}=\id_{E}$. 
One can expect that there is a semilocal analog of the Conn
linearization theorem in the case when the transverse Lie
algebra $\mathfrak{g}$ of the leaf~$B$ is semisimple of compact
type. 
But, according to an important observation due to~\cite{DW},
there are \textit{nonlinearizable} Poisson structures even if
$\mathfrak{g}$ is semisimple of compact type. The corresponding
examples are derived in the class of Poisson structures called
\textit{Casimir-weighted products}~\cite{DW}. 
This result rises the question
on describing the corresponding obstructions to the
linearizability. Here we make an attempt to answer this
question. 

Let $\Casim_{B}(E)$ be the space of Casimir functions of the
linearized transverse Poisson structure $\Lambda$ vanishing on~$B$.

\begin{claim} 
There exists an intrinsic coboundary operator
$$
\partial_{0}:\Omega^{k}(B)\otimes\Casim_{B}(E)\rightarrow
\Omega^{k+1}(B)\otimes\Casim_{B}(E)
$$
associated to $(\Psi,B)$.
\end{claim}

One can define $\partial_{0}$ in terms of the exterior covariant
derivative $\partial^{\Gamma^{(1)}}$ of a homogeneous connection
$\Gamma^{(1)}$, but this definition is independent of
$\Gamma^{(1)}$ and the operator $\partial_{0}$ is rather
attributed to the (singular) symplectic foliation of the
linearized transverse Poisson structure $\Lambda$. Next, we can
associated to the pair $(\Pi,\Pi^{(1)})$ a $\partial_{0}$-cocycle
$C=C_{\Pi,\Pi^{(1)}}\in\Omega^{k}(B)\otimes\Casim_{B}(E)$,
called a \textit{Casimir $2$-cocycle}. The $2$-form
$C_{\Pi,\Pi^{(1)}}$ is computed in terms of the geometric data
$(\Pi_{V},\Gamma,F)$ and $(\Lambda,\Gamma^{(1)},F^{(1)})$ and
is well defined under the assumption on the triviality of the
first reduced Poisson cohomology space of $\Lambda$. In
particular, 
this assumption holds in the case when the transverse Lie
algebra $\mathfrak{g}$ is semisimple of compact type. We show that
the $\partial_{0}$-cohomology class $[C_{\Pi,\Pi^{(1)}}]$ is an
invariant of the leaf $B$. Moreover, using a Casimir $2$-cocycle
$C=C_{\Pi,\Pi^{(1)}}$, we can define the \textit{deformed
coupling form} $^{(1)}F+C_{\Pi,\Pi^{(1)}}$ and the
\textit{deformed linearized Poisson structure} $\Pi_{C}^{(1)}$
as the coupling Poisson tensor associated with the data
$(\Lambda,\Gamma^{(1)},F^{(1)}+C_{\Pi,\Pi^{(1)}})$. Assuming
that the transverse Lie algebra $\mathfrak{g}$ of the leaf $B$ 
is \textit{semisimple of compact type}, we formulate our main
results as follows.

\begin{claim}
{\rm(}Normal Form Theorem{\rm)}. 
The germs at $B$ of the Poisson structures $\Psi$ and $\Pi_{C}^{(1)}$
are isomorphic. 
\end{claim}

A similar result to this statement was formulated in~\cite{Br} 
but without determining the deformed coupling form.

\begin{claim}
{\rm(}Semilocal Linearization Theorem{\rm)}. 
The Poisson structure $\Psi$ is linearizable at $B$ if and only
if the $\partial_{0}$-cohomology class of the leaf $B$ is zero,
$$
[C_{\Pi,\Pi^{(1)}}]=0.
$$
\end{claim}

These results are based on the
following observation~\cite{Vo1}: a homotopy for coupling
Poisson structures implies the Poisson equivalence. 
Here we give a further development of this thesis for general
families of Poisson tensors. 

We also refer to the work~\cite{CrFe}, where one can find a
conjecture on the linearization around a compact symplectic leaf
$B$, formulated in terms of the integrability of the transitive
Lie algebroid of~$B$. 

The paper is organized as follows.

In Sections 2.1--2.3 we give main definitions and formulate some
useful technical results in Propositions~2.11 and~2.13. 
In Section~2.4, we study the relationship between a general
homotopy for coupling structures and the Poisson equivalence. 
Here the main results are formulated in Theorems~2.18 and~2.19. 
In Section~2.5, we introduce the notion of a relative Casimir 
$2$-cocycle of two coupling Poisson tensors and discuss the
corresponding properties. 
In Section~3, using results of Section~2, we derive some
criteria for the equivalence of Poisson structures near a given
symplectic leaf. First, in Theorem~3.2, we show that the
equivalence of coupling Poisson tensors imply the equivalence of
the corresponding vertical parts. Sufficient conditions for the 
Poisson equivalence are presented in Theorems~3.4 and~3.8 
and Proposition~3.9. We show also (Theorem~3.11) that the
cohomology class of the Casimir $2$-cocycle can be viewed as
obstruction to the equivalence of coupling Poisson tensors. 
Section~4 is devoted to the semilocal linearization problem.
The main results are presented in Theorems~4.12 and~4.14.

\textbf{Acknowledgement.}
I am grateful M.~V.~Karasev for helpful and stimulating
discussions during the preparation of the text. 

\section{Poisson coupling}

In this section we first recall some basic properties of
coupling Poisson structures on fibered spaces~\cite{Vo1}. For
generalizations of the coupling procedure to foliated manifolds
and Jacobi manifolds, we refer to~\cite{Va2}. We introduced a
homotopy equivalence for coupling Poisson tensors and formulate
some useful technical results. Finally, we introduce the notion
of the Casimir $2$-cocycle of two coupling Poisson tensors 
which will play an important role in formulating the main results
in Sections~3 and~4.

We shall perform all computations in coordinates. A
free-coordinate approach can be found in~\cite{Va2}.

\subsection{Preliminaries}
Let us fix some notation and recall the basic definitions. 
Let $\pi:\,E\rightarrow B$ be a fiber bundle over a base~$B$. 
By $\chi^{k}(E)$ and $\Omega^{k}(E)$ we
will denote the space of antisymmetric $k$-tensor fields and
$k$-forms on the total space $E$. In particular, 
$\mathcal{X}(E)=\chi^{1}(E)$ is the space of
smooth vector fields on $E$. Let $\mathbb{V}=\ker d\pi\subset
TE$ be the vertical subbundle. The elements of $\chi^{k}(E)$,
tangent to $\mathbb{V}$, 
form the subspace of vertical $k$-tensor fields denoted by
$\chi_{V}^{k}(E)=\Sect(\wedge^{k}\mathbb{V})$. 
A $k$-form on $E$ is said to be horizontal if it annihilates 
the vertical subbundle. We denote the
subspace of horizontal $k$-forms by $\Omega_{H}^{k}(E)$.

Suppose we are given an \textit{Ehresmann connection}~\cite{GHV}
on $E$, that is, a smooth splitting
\begin{equation}
TE=\mathbb{H}\oplus\mathbb{V}
\tag{2.1}
\end{equation}
which is given by a subbundle $\mathbb{H}\subset TE$ 
called \textit{horizontal}.
The $k$-vector fields tangent to $\mathbb{H}$ are called
horizontal and form the subspace 
$\chi_{H}^{k}(E)\approx\Sect(\wedge^{k}\mathbb{H})$. 
In particular, $\mathcal{X}_{H}(E)=\chi_{H}^{1}(E)$ is the space
of horizontal vector fields. For every vector field $X$ on $E$
we have the decomposition $X=X_{H}+X_{V}$ onto the horizontal
and vertical components $X_{H}$ and $X_{V}$, respectively. 
More generally,
$$ 
\chi^{k}(E)=\bigoplus_{i+j=k}\chi_{H}^{i}(E)\wedge\chi_{V}^{j}(E),
$$
where an element of $\chi_{H}^{i}(E)\wedge\chi_{V}^{j}(E)$ is
said to be a multivector field of degree~$(i,j)$.

Let $\Gamma\in\Omega^{1}(E,\mathbb{V})$ be the connection form
of (2.1), $\Gamma(X)\overset{\text{def}}{=}X_{V}$. 
Thus, $\mathbb{H}=\ker\Gamma$. 
The curvature form 
$\Curv^{\Gamma}\in\Omega^{2}(B,\End(\mathbb{V}))$ is given by
$$ 
\Curv^{\Gamma}(u_{1},u_{2})=\hor^{\Gamma}
([u_{1},u_{2}])-[\hor^{\Gamma}(u_{1}),\hor^{\Gamma}(u_{2})].
$$
Here $\hor^{\Gamma}(u)$ is the $\Gamma$-\textit{horizontal lift}
of a vector field $u\in\mathcal{X}(B)$.

Consider the tensor product $\Omega^{k}(B)\otimes C^{\infty}(E)$ over
$C^{\infty}(B)$. The $\Gamma$-\textit{co\-var\-iant exterior derivative\/}
$\partial^{\Gamma}:\,\Omega^{k}(B)\otimes C^{\infty}(E)
\rightarrow\Omega^{k+1}(B)\otimes C^{\infty}(E)$ is defined as
\begin{align} 
& (\partial^{\Gamma}\Theta)(u_{0},u_{1},\dots,u_{k})
\overset{\text{def}}{=}
\sum_{i=0}^{k}(-1)^{i} L_{\hor(u_{i})}
\Theta(u_{0},u_{1},\dots,\hat{u}_{i},\dots,u_{k})
\tag{2.2}\\
&\qquad 
+\sum_{0\leq i<j\leq k}(-1)^{i+j}
\Theta\big([u_{i},u_{j}],u_{0},u_{1},\dots,\hat{u}_{i},\dots,
\hat{u}_{j},\dots,u_{k}\big).
\nonumber
\end{align}
We have the identity 
\begin{align}
& \big((\partial^{\Gamma})^{2}\Theta)
(u_{0},u_{1},\dots,u_{i},\dots,u_{j},\dots,u_{k+1}\big)
\tag{2.3}\\
&\qquad 
=\sum_{0\leq i<j\leq k+1}(-1)^{i+j}
L_{\Curv^{\Gamma}(u_{i},u_{j})}
\Theta(u_{0},u_{1},\dots,\hat{u}_{i},\dots,\hat{u}_{j},\dots,u_{k+1})
\nonumber
\end{align}
which says that the coboundary condition
$(\partial^{\Gamma})^{2}=0$ holds if and only if 
$\Curv^{\Gamma}=0$. This is the integrability
condition for the horizontal plane distribution.

One can assign to every 
$\Theta\in\Omega^{k}(B)\otimes C^{\infty}(E)$ 
a horizontal $k$-form $\pi^{\ast}\Theta$ 
uniquely determined by the condition
$$
(\pi^{\ast}\Theta)
\big(\hor^{\Gamma}(u_{1}),\dots,\hor^{\Gamma}(u_{k})\big)
=\Theta(u_{1},\dots,u_{k})
$$
for any $u_{1},\dots,u_{k}\in\mathcal{X}(B)$. 
This correspondence is independent of the choice of a connection.

Consider a (local) coordinate system
$(\xi,x)=(\xi^{i},x^{\sigma})$ on the total space $E$, 
where $(\xi^{i})$ are coordinates on the base $B$ and
$(x^{\sigma})$ are coordinates along the fibers of~$E$.  
In coordinates, we have
$$
\Gamma=\Gamma^{\nu}\otimes\frac{\partial}{\partial x^{\nu}},
\qquad 
\Gamma^{\nu}=dx^{\nu}+\Gamma_{i}^{\nu}(\xi,x)d\xi^{i}.
$$
Here and throughout the text, the summation over repeated
indices will be understood. Locally, the horizontal subbundle is
generated by the vector fields 
$$
\hor_{i}\overset{\text{def}}{=}\hor^{\Gamma}
\bigg(\frac{\partial}{\partial\xi^{i}}\bigg)
=\frac{\partial}{\partial\xi^{i}}-\Gamma_{i}^{\nu}(\xi,x)
\frac{\partial}{\partial x^{\nu}}.
$$
Moreover, the curvature form can be viewed as a $2$-form on $B$
with values in the space of vertical vector fields,
$$
\Curv^{\Gamma}=\frac{1}{2}\Curv_{ij}^{\sigma}
(\xi,x)d\xi^{i}\wedge d\xi^{j}\otimes
\frac{\partial}{\partial x^{\sigma}}.
$$
Notice that
\begin{equation}
[\hor_{i},\hor_{j}]=-\Curv_{ij}^{\sigma}
\frac{\partial}{\partial x^{\sigma}}\in\mathcal{X}_{V}(E),
\tag{2.4}
\end{equation}
and
\begin{equation}
\bigg[\hor_{i},\frac{\partial}{\partial x^{\sigma}}\bigg]
=\frac{\partial\Gamma_{i}^{\nu}}{\partial x^{\sigma}}
\frac{\partial}{\partial x^{\nu}}\in\mathcal{X}_{V}(E).
\tag{2.5}
\end{equation}

Recall some properties of the \textit{Schouten bracket} for
multivector fields on $E$ (see, for example \cite{KM,Va1}). 
Let $\Pi\in\chi^{2}(E)$ be a bivector
field and $\Pi^{\natural}:\,T^{\ast}E\rightarrow TE$
be the induced bundle morphism, 
$\langle\beta,\Pi^{\natural}(\alpha)\rangle=\Pi(\alpha,\beta)$. 
For any $G\in C^{\infty}(E)$ and $W\in\mathcal{X}(E)$, 
we have $[G,\Phi]=\Pi^{\natural}dG$
and $[W,\Pi]=L_{W}\Pi$, where $L_{W}$ is the Lie derivative,
\begin{equation}
L_{W}\Pi(\alpha,\beta)
=L_{\Pi^{\natural}(\alpha)}\beta(W)
-L_{\Pi^{\natural}(\beta)}\alpha(W)+L_{W}(\Pi(\alpha,\beta)).
\tag{2.6}
\end{equation}
Moreover, we also need the following identities:
\begin{equation}
[ W,Z_{1}\wedge Z_{2}]
=[W,Z_{2}]\wedge Z_{1}-[W,Z_{1}]\wedge Z_{2}
\tag{2.7}
\end{equation}
and
\begin{equation}
[Z_{1}\wedge Z_{2},\Pi]=Z_{2}\wedge L_{Z_{1}}\Pi-Z_{1}\wedge 
L_{Z_{2}}\Pi
\tag{2.8}
\end{equation}
for $W,Z_{1},Z_{2}\in\mathcal{X}(E)$.

\textbf{The Jacobi identity}. 
Fix an Ehresmann connection $\Gamma$. Suppose we
are given a bivector field $\Pi\in\chi^{2}(E)$ which has the following
representation with respect to decomposition~(2.1):
\begin{equation}
\Pi=\Pi_{H}+\Pi_{V}.
\tag{2.9}
\end{equation}
Here $\Pi_{H}\in\chi_{H}^{2}(E)$ and $\Pi_{V}\in\chi_{V}^{2}(E)$ 
are horizontal and vertical bivector fields, respectively,
\begin{equation}
\Pi_{H}=\frac{1}{2}\Pi_{H}^{ij}(\xi,x)\hor_{i}\wedge\hor_{j},
\qquad 
\Pi_{V}=\frac{1}{2}\Pi_{V}^{\alpha\beta}(\xi,x)
\frac{\partial}{\partial x^{\alpha}}\wedge
\frac{\partial}{\partial x^{\beta}}.
\tag{2.10}
\end{equation}
Let $\mathbb{V}^{0}\overset{\text{def}}{=}\Ann(\mathbb{V})\subset
T^{\ast}E$ be the annihilator of the vertical subbundle
$\mathbb{V}$. Notice that (2.9) holds if for every horizontal
$1$-form $\beta$ the image $\Pi^{\natural}(\beta)$ 
is a horizontal vector field,
\begin{equation}
\Pi^{\natural}(\mathbb{V}^{0})\subseteq\mathbb{H}.
\tag{2.11}
\end{equation}
This is just a condition for vanishing the component of degree
$(1,1)$ in the decomposition of~$\Pi$.

\begin{proposition}
The Jacobi identity for the bivector field $\Pi$,
\begin{equation}
[\Pi,\Pi]=0
\tag{2.12}
\end{equation}
is equivalent to the following relations:

\begin{align}
\underset{(\alpha,\beta,\gamma)}{\mathfrak{S}}
\Pi_{V}^{\alpha\sigma}
\frac{\partial\Pi_{V}^{\beta\gamma}}{\partial x^{\sigma}}&=0,
\tag{2.13}
\\
\Pi_{H}^{is}L_{\hor_{s}}\Pi_{V}&=0,
\tag{2.14}
\\
\underset{(k,i,j)}{\mathfrak{S}}
\Pi_{H}^{ks} L_{\hor_{s}}(\Pi_{H}^{ij})&=0,
\tag{2.15}
\\
\Pi_{H}^{is}\Curv_{sj^{\prime}}^{\sigma}\Pi_{H}^{j^{\prime}j}
&=\Pi_{V}^{\sigma\sigma^{\prime}}
\frac{\partial}{\partial x^{\sigma^{\prime}}} 
(\Pi_{H}^{ij}).
\tag{2.16}
\end{align}

Here $\mathfrak{S}$ denotes the cyclic sum and $L_{\hor_{s}}$ is the
Lie derivative along $\hor_{s}$.
\end{proposition}

\begin{proof}
Using relations (2.4)--(2.8), we compute the Schouten bracket
$[\Pi,\Pi]$ by parts. 

(i) The ``horizontal--horizontal'' part:
\begin{align*}
[\Pi_{H},\Pi_{H}]  
&=\frac{1}{4}[\Pi_{H}^{ij}\hor_{i}\wedge\hor_{j},
\Pi_{H}^{i^{\prime}j^{\prime}}\hor_{i^{\prime}}\wedge\hor_{j^{\prime}}]
\\
&\qquad
-\Pi_{H}^{ks}L_{\hor_{s}}(\Pi_{H}^{ij})\hor_{i}\wedge\hor_{j}\wedge\hor_{k}
\\
&\qquad
-\Pi_{H}^{ji}\Curv_{^{ii^{\prime}}}^{\sigma}
\Pi_{H}^{i^{\prime}j^{\prime}}\frac{\partial}{\partial x^{\sigma}}
\wedge\hor_{j}\wedge\hor_{j^{\prime}}.
\end{align*}

(ii) The ``horizontal--vertical'' part:
\begin{align*}
2[\Pi_{H},\Pi_{V}]  
&=\frac{1}{2}
\bigg[\Pi_{H}^{ij}\hor_{i}\wedge\hor_{j},
\Pi_{V}^{\sigma\sigma^{\prime}}\frac{\partial}{\partial x^{\sigma}}
\wedge\frac{\partial}{\partial x^{\sigma^{\prime}}}\bigg]
\\
&=\Pi_{V}^{\sigma^{\prime}\sigma}
\frac{\partial\Pi_{H}^{jj^{\prime}}}{\partial x^{\sigma^{\prime}}}
\frac{\partial}{\partial x^{\sigma}}
\wedge\hor_{j}\wedge\hor_{j^{\prime}}
-\frac{1}{2}\Pi_{H}^{ij}\hor_{i}\wedge L_{\hor_{j}}\Pi_{V}.
\end{align*}
Here
$$
L_{\hor_{j}}\Pi_{V}=\frac{1}{2}
\bigg(L_{\hor_{j}}\Pi_{V}^{\alpha\beta}
+\Pi_{V}^{\alpha\sigma}
\frac{\partial\Gamma_{j}^{\beta}}{\partial x^{\sigma}}
-\Pi_{V}^{\beta\sigma}
\frac{\partial\Gamma_{j}^{\alpha}}{\partial x^{\sigma}}\bigg)
\frac{\partial}{\partial x^{\alpha}}\wedge
\frac{\partial}{\partial x^{\beta}}.
$$

(iii) The ``vertical--vertical'' part:
\begin{align*}
[\Pi_{V},\Pi_{V}]  
&=\frac{1}{4}\bigg[\Pi_{V}^{\sigma\sigma^{\prime}}
\frac{\partial}{\partial x^{\sigma}}\wedge
\frac{\partial}{\partial x^{\sigma^{\prime}}},
\Pi_{V}^{\alpha\beta}\frac{\partial}{\partial x^{\alpha}}\wedge
\frac{\partial}{\partial x^{\beta}}\bigg]
\\
&=\Pi_{V}^{\sigma\nu}
\frac{\partial\Pi_{V}^{\alpha\beta}}{\partial x^{\nu}}
\frac{\partial}{\partial x^{\sigma}}\wedge
\frac{\partial}{\partial x^{\alpha}}\wedge
\frac{\partial}{\partial x^{\beta}}.
\end{align*}
It follows from here that the Schouten bracket $[\Pi,\Pi]$ is
the sum of $3$-vector fields of degrees $(3,0)$, $(2,1)$, $(1,2)$, 
and $(0,3)$. Vanishing of these terms leads to relations 
(2.13)--(2.16).
\end{proof}

\begin{corollary}
The Jacobi identity for $\Pi$ in {\rm(2.9)} implies that its
vertical part $\Pi_{V}$ is a Poisson tensor.
\end{corollary}

\begin{remark}
In \cite{Va2}, a Poisson tensor satisfying condition (2.11) is
called \textit{almost coupling}.
\end{remark}

\subsection{Coupling Poisson tensors}
A bivector field $\Pi\in\chi^{2}(E)$ is said to be
\textit{horizontally nondegenerate} on $E$ if for every $m\in E$
the image $\Pi_{H}^{\natural}(\mathbb{V}_{m}^{0})$ is a
complementary subspace of $\mathbb{V}_{m}$ in $T_{m}E$. 
In this case, $\Pi$ induces the Ehresmann
connection $\Gamma$ corresponding to the horizontal subbundle
\begin{equation}
\mathbb{H}\overset{\text{def}}{=}\Pi^{\natural}(\mathbb{V}^{0}).
\tag{2.17}
\end{equation}
This implies that $\Pi(\alpha,\beta)=0$ for all 
$\alpha\in\Sect(\mathbb{H}^{0})$ and $\beta\in\Sect(\mathbb{V}^{0})$.
Consequently, $\Pi$ has representation (2.9) with respect to
connection (2.17). Then, the horizontal part $\Pi_{H}$ is
nondegenerate in the sense that
\begin{equation}
\ker\Pi_{H}^{\natural}|_{\mathbb{V}_{m}^{0}}=\{0\},
\tag{2.18}
\end{equation}
or, equivalently, 
$\Pi_{H}^{\natural}(\mathbb{V}_{m}^{0})=\mathbb{H}_{m}$. 
It follows from here that there exists a unique $2$-form
$F\in\Omega^{2}(B)\otimes C^{\infty}(E)$ defined by the condition
\begin{equation}
\Pi_{H}(\beta_{1},\beta_{2})
=(\pi^{\ast}F)(\Pi_{H}^{\natural}(\beta_{1}),\Pi_{H}^{\natural}(\beta_{2}))
\tag{2.19}
\end{equation}
for any $\beta_{1},\beta_{2}\in\Sect(\mathbb{V}^{0})$. 
It is clear that
\begin{equation}
\pi^{\ast}F |_{\mathbb{H}_{m}}\quad\text{\textit{is nondegenerate}}
\tag{2.20}
\end{equation}
for every $m\in E$. In coordinates,
$$
\pi^{\ast}F=F_{ij}(\xi,x)d\xi^{i}\wedge d\xi^{j},
$$
where $\Pi_{H}^{is}(\xi,x) F_{sj}(\xi,x)=-\delta_{j}^{i}$ and
\begin{equation}
\det(F_{ij}(\xi,x))\neq0.
\tag{2.21}
\end{equation}
Thus, the horizontally nondegenerate bivector field $\Pi$
induces the triple $(\Pi_{V},\Gamma,F)$ consisting of the
vertical bivector field $\Pi_{V}$, the Ehresmann connection
$\Gamma$ in (2.17) and the $2$-form $F$ defined by (2.19). 
The triple $(\Pi_{V},\Gamma,F)$ will be called the
\textit{geometric data} of~$\Pi$. This correspondence is
one-to-one. 
Conversely, for a given $(\Pi_{V},\Gamma,F)$, where $F$
satisfies the nondegeneracy condition (2.20), the corresponding
horizontally nondegenerate bivector field $\Pi$ is defined by
\begin{equation}
\Pi=-\frac{1}{2}F^{ij}(\xi,x)\hor_{i}\wedge\hor_{j}
+\frac{1}{2}\Pi_{V}^{\alpha\beta}(\xi,x)
\frac{\partial}{\partial x^{\alpha}}\wedge
\frac{\partial}{\partial x^{\beta}}.
\tag{2.22}
\end{equation}
Here, $F^{is}(\xi,x)F_{sj}(\xi,x)=\delta_{j}^{i}$.

Taking into account (2.21), we deduce the following result
straightforwardly from Proposition~2.1.

\begin{proposition}
A horizontally nondegenerate bivector field $\Pi$ is a Poisson
tensor if and only if the geometric data $(\Pi_{V},\Gamma,F)$
satisfy the following conditions
\begin{align}
[\Pi_{V},\Pi_{V}]&=0,
\tag{2.23}
\\
L_{\hor^{\Gamma}(v)}\Pi_{V}&=0,
\tag{2.24}
\\
\partial^{\Gamma}F&=0,
\tag{2.25}
\\
\Curv^{\Gamma}(v,u)&=(\Pi_{V})^{\sharp}d(F(v,u))
\tag{2.26}
\end{align}
for any $v,u\in\mathcal{X}(B)$.
\end{proposition}

By (2.16) and (2.21) we also get  the following fact.

\begin{corollary}
The Jacobi identity for the horizontal component $\Pi_{H}$ of $\Pi$
is equivalent to the zero curvature condition,
\begin{equation}
\Curv^{\Gamma}=0.
\tag{2.27}
\end{equation}
\end{corollary}

The horizontally nondegenerate Poisson tensor $\Pi$ in (2.22) is
called a \textit{coupling Poisson structure} associated with
geometric data $(\Pi_{V},\Gamma,F)$.  
The $2$-form $F$ will also be called a \textit{coupling form}.
Formula (2.22) is the result of coupling $F$ and $\Pi_{V}$ via
the connection~$\Gamma$.

The pairwise Poisson brackets of coordinate functions relative
to $\Pi$ are written in terms of geometric data as follows:
\begin{align}
\{\xi^{i},\xi^{j}\}_{\Pi}&=-F^{ij}(\xi,x),
\tag{2.28}
\\
\{\xi^{i},x^{\sigma}\}_{\Pi}&=F^{is}(\xi,x)\Gamma_{s}^{\sigma}(\xi,x),
\tag{2.29}
\\
\{x^{\alpha},x^{\beta}\}_{\Pi}&=\Pi_{V}^{\alpha\beta}(\xi,x)
-\Gamma_{i}^{\alpha}(\xi,x)F^{ij}(\xi,x)\Gamma_{j}^{\beta}(\xi,x).
\tag{2.30}
\end{align}

\begin{example}
Let $E=B\times N$ is the product of a symplectic manifold
$(B,\omega)$ and a Poisson manifold $N$ with Poisson tensor
$\Phi$. Let $p_{1}$ and $p_{2}$ be the canonical projections to
the first and the second factors, respectively. 
We can think of $E$ as a trivial fiber bundle over $B$ with 
$\pi=p_{1}$.
Then the direct product Poisson structure on $B\times N$ is a
coupling Poisson structure associated with the following data:
$(\Pi_{V}$, $\mathbb{H}=\ker p_{2}$, $F=\omega\otimes1)$. 
Here $\Pi_{V}$ is a vertical bivector field on $B\times N$ which
is uniquely determined by the condition: 
$\Pi_{V}(p_{2}^{\ast}\alpha,p_{2}^{\ast}\beta)
=p_{2}^{\ast}(\Phi(\alpha,\beta))$ for every 
$\alpha,\beta\in\Omega^{1}(N)$. In this case, $\Gamma$ is a flat
connection with trivial holonomy. By the local splitting
theorem~\cite{We1}, locally, each Poisson manifold has this form.
\end{example}

\begin{example}
Every Poisson manifold in a neighborhood of a symplectic leaf is
realized as a coupling Poisson structure on the normal bundle
(see Section~4). 
\end{example}

More examples can be found in \cite{Vo1,Vo4,Va2,Va3,DW}.

\textbf{Fiber-tangent Poisson structures.}
By (2.23) the vertical part of a coupling Poisson tensor is also
Poisson. Following the terminology introduced in~\cite{Va2}, 
a vertical Poisson tensor $\Upsilon$ on the total space $E$ will
be called a \textit{fiber-tangent Poisson structure}. 
Such a Poisson structure can be uniquely characterized by the property:  
for ever $m\in E$ the symplectic leaf of $\Upsilon$ through $m$ 
belongs to the fiber $E_{\pi(m)}$. 
The restriction $\Upsilon_{\xi}=\Upsilon |_{E_{\xi}}$
is well defined and gives a Poisson tensor on the fiber
$E_{\xi}$ which varies smoothly with $\xi\in B$. 
Thus, one can think of $\pi:\,E\rightarrow B$ as a bundle of
Poisson manifolds with fiberwise Poisson structure $\Upsilon_{\xi}$. 
If $\Upsilon$ is the vertical part of a coupling Poisson tensor
$\Pi$, $\Upsilon=\Pi_{V}$, then condition (2.24) implies that the
parallel transport operator of the connection $\Gamma$ preserves
the fiberwise Poisson structure on $E$ induced by $\Upsilon$. 
Such a connection is called a \textit{Poisson connection}.

\begin{definition}
A fiber-tangent Poisson structure $\Upsilon$ is said to be
\textit{locally trivial} if in a neighborhood of each point in
$E$ there exists a coordinate system $(\xi,x)=(\xi^{i},x^{\sigma})$ 
such that $\Upsilon$ has the form
\begin{equation}
\Upsilon=\frac{1}{2}\Upsilon^{\alpha\beta}(x)
\frac{\partial}{\partial x^{\alpha}}\wedge
\frac{\partial}{\partial x^{\beta}}.
\tag{2.31}
\end{equation}
Here $(\xi^{i})$ are coordinates on the base $B$ and
$(x^{\sigma})$ are coordinates on the fibers of~$E$.
\end{definition}

In the next section, we show that a typical feature of the
fiber-tangent Poisson structure $\Upsilon$ which comes from a
coupling Poisson tensor is  that $\Upsilon$ is locally trivial.
The local triviality property allows us to transfer some
properties of $\Upsilon$ from the fiber to the total space using
the partition of unity argument. 

\textbf{Fiber preserving transformations}. 
Suppose we are given a coupling tensor $\Pi$ on $E$ associated
with geometric data $(\Pi_{V},\Gamma,F)$. 

\begin{proposition}
Let $g$ be a fiber preserving diffeomorphism. Then the push-forward
$\tilde{\Pi}=g_{\ast}\Pi$ is a coupling Poisson tensor
associated to the geometric data
\begin{equation}
\tilde{\Pi}_{V}=g_{\ast}\Pi_{V},\qquad
\tilde{\Gamma}=g_{\ast}\Gamma,\qquad
\tilde{F}=g_{\ast}F.
\tag{2.32}
\end{equation}
\end{proposition}

\begin{proof} 
Consider the bivector field $\tilde{\Pi}=g_{\ast}\Pi$ and the
corresponding distribution $\tilde{\mathbb{H}}$ defined by (2.17),
\begin{equation}
\tilde{\mathbb{H}}_{g(m)}
=(d_{m}g)\circ\Pi_{m}^{\sharp}\circ(d_{m}g)^{\ast}(\mathbb{V}_{g(m)}^{0})
\tag{2.33}
\end{equation}
for every $m\in N$. Here $d_{m}g:\,T_{m}E\rightarrow T_{g(m)}E$
is the tangent map of $g$. On the other hand, consider the
push-forward $g_{\ast}\mathbb{H}$ of $\mathbb{H}$ by $g$,
\begin{equation}
(g_{\ast}\mathbb{H)}_{g(m)}=(d_{m}g)\mathbb{H}_{m}.
\tag{2.34}
\end{equation}
By assumption, $g$ is fiber preserving, that is, $dg$ preserves
the vertical subbundle,
$$
(d_{m}g)\mathbb{V}_{m}=\mathbb{V}_{g(m)}
$$
and hence 
$(d_{m}g)^{\ast}(\mathbb{V}_{g(m)}^{0}=\mathbb{V}_{m}^{0}$.
Comparing (2.33) and (2.34) leads to
$$
\tilde{\mathbb{H}}=g_{\ast}\mathbb{H}.
$$
Thus, $\tilde{\mathbb{H}}$ is complementary to $\mathbb{V}$ and
hence $\tilde{\Pi}$ is a coupling Poisson tensor. For the
connection form we have 
$\tilde{\Gamma}((d_{m}g)X)=(d_{m}g)\Gamma(X)$ 
for $X\in T_{m}E$. Other two identities in (2.32) are evident.
\end{proof}

It follows that decomposition (2.9) is stable under a fiber
preserving diffeomorphism $g$,
\begin{equation}
g_{\ast}(\Pi_{H}+\Pi_{V})=(g_{\ast}\Pi)_{H}+(g_{\ast}\Pi)_{V}.
\tag{2.35}
\end{equation}
In general, this is not true.

\subsection{Infinitesimal Poisson automorphisms}
Let $\Pi$ be a coupling Poisson tensor associated with geometric
data $(\Pi_{V},\Gamma,F)$. Our goal is to describe infinitesimal
Poisson automorphisms (Poisson vector fields) of $\Pi$ in terms
of the geometric data.

We start with some useful technical formulas. 
Let
\begin{equation}
X=X^{i}(\xi,x)\hor_{i}
\tag{2.36}
\end{equation}
be a horizontal vector field. Denote 
$\mathbf{i}_{X}F\in\Omega^{1}(B)\otimes C^{\infty}(N)$ 
defined by 
$\pi^{\ast}(\mathbf{i}_{X}F)=\mathbf{i}_{X}(\pi^{\ast}F)$, 
or locally, 
$\mathbf{i}_{X}F=X^{s}(\xi,x)F_{sj}(\xi,x)\otimes d\xi^{j}$.

\begin{lemma}
For any $v_{1},v_{2}\in\mathcal{X}(B)$, we have
\begin{equation}
L_{X}(\pi^{\ast}F)\big(\hor(v_{1}),\hor(v_{2})\big)
=\partial^{\Gamma}(\mathbf{i}_{X}F)(v_{1},v_{2}).
\tag{2.37}
\end{equation}
\end{lemma}

\begin{proof}
In components, condition (2.25) for $F$ reads
$$
\underset{(ijs)}{\mathfrak{S}}L_{\hor_{i}}F_{js}=0.
$$
Applying this identity and using the definition of 
$\partial^{\Gamma}$, we get
\begin{align*}
L_{X}(F)(\hor_{m},\hor_{m^{\prime}})  
&=L_{X}F_{mm^{\prime}}-F_{sm}L_{\hor_{m^{\prime}}}X^{s}
+F_{sm^{\prime}}L_{\hor_{m}}X^{s}
\\
&=L_{\hor_{m}}(\mathbf{i}_{X}F)_{m^{\prime}}
-L_{\hor_{m^{\prime}}}(\mathbf{i}_{X}F)_{m}
\\
&=(\partial^{\Gamma}(\mathbf{i}_{X}F))_{mm^{\prime}}.
\end{align*}
\end{proof}

\begin{proposition}
The Lie derivative of the coupling Poisson tensor $\Pi$ along a
horizontal vector field $X$ is given by the formula
\begin{align}
L_{X}\Pi &=\frac{1}{2}
F^{im}[(\partial^{\Gamma}(\mathbf{i}_{X}F))_{mm^{\prime}}]
F^{m^{\prime}j}\hor_{i}\wedge\hor_{j}
\tag{2.38}
\\
&\qquad 
-F^{si}\hor_{i}\wedge\Pi_{V}^{\sharp}(d(X^{m}F_{ms})).
\nonumber
\end{align}
\end{proposition}

\begin{proof}
It follows from (2.4) and (2.26) that
$$
[\hor_{s},\hor_{m}]
=-\Curv_{sm}^{\sigma}\frac{\partial}{\partial x^{\sigma}}
=-\Pi_{V}^{\sharp}(dF_{sm}).
$$
Then,
$$
L_{\hor_{s}}\Pi_{H}=-\frac{1}{2}(L_{\hor_{s}}F^{ij})\hor_{i}
\wedge\hor_{j}-F^{mi}\hor_{i}\wedge\Pi_{V}^{\sharp}(dF_{sm})
$$
and
\begin{align*}
L_{X}\Pi_{H}  
& =\frac{1}{2}F^{im}(L_{X}F_{mm^{\prime}}-F_{sm}
L_{\hor_{m^{\prime}}}X^{s}+F_{sm^{\prime}}L_{\hor_{m}}X^{s})
F^{m^{\prime}j}\hor_{i}\wedge\hor_{j}
\\
&\qquad 
-X^{s}F^{mi}\hor_{i}\wedge\Pi_{V}^{\sharp}(dF_{sm})
-\hor_{i}\wedge(\Pi_{V}^{\sharp}dX^{i})
\end{align*}
From here, by using (2.37), we deduce (2.38).
\end{proof}

Denote by $\Casim(E,\Pi_{V})$ the space of Casimir functions of
the vertical Poisson structure $\Pi_{V}$. It is clear that 
$\pi^{\ast}C^{\infty}(B)\subset\Casim(E,\Pi_{V})$.

\begin{corollary}
A horizontal vector $X$ is an infinitesimal Poisson automorphism
of $\Pi$ if and only if
\begin{equation}
L_{X}(\pi^{\ast}F) |_{\mathbb{H}}=0,
\tag{2.39}
\end{equation}
and
\begin{equation}
\mathbf{i}_{X}F\in\Omega^{1}(B)\otimes\Casim(E,\Pi_{V}).
\tag{2.40}
\end{equation}
\end{corollary}

Condition (2.39) says that the Lie derivative of the horizontal
$2$-form $\pi^{\ast}F$ along $X$ is a vertical form.

Furthermore, using (2.24), we get the formula
\begin{equation}
L_{X}\Pi_{V}=-\hor_{i}\wedge(\Pi_{V}^{\sharp}dX^{i})
\tag{2.41}
\end{equation}
which says that $X$ is an infinitesimal Poisson automorphism of
$\Pi_{V}$ if and only if the coefficients $X^{i}$ in (2.36) are
(local) Casimir functions. 

Now, let
$$
Y=Y^{\sigma}(\xi,x)\frac{\partial}{\partial x^{\sigma}}
$$
be a vertical vector field. By straightforward calculation we get the
formula
\begin{align}
L_{Y}\Pi 
&=\frac{1}{2} F^{im}[L_{Y}F_{mm^{\prime}}]F^{m^{\prime}j}
\hor_{i}\wedge\hor_{j}
\tag{2.42}
\\
&\qquad 
-F^{si}[\hor_{s},Y]^{\sigma}\hor_{i}\wedge
\frac{\partial}{\partial x^{\sigma}}+L_{Y}\Pi_{V},
\nonumber
\end{align}
which gives the following criterion.

\begin{proposition}
A vertical vector field $Y$ is an infinitesimal Poisson
automorphism of $\Pi$ if and only if
\begin{align}
L_{Y}\Pi_{V}&=0,
\tag{2.43}
\\
L_{Y}(\pi^{\ast}F)&=0,
\tag{2.44}
\\
[\hor(u),Y]&=0
\tag{2.45}
\end{align}
for every $u\in\mathcal{X}(B)$.
\end{proposition}

Remark that condition (2.45) means that the horizontal subbundle
$\mathbb{H}$ is invariant under the flow of~$Y$.

\subsection{Poisson homotopy}
Suppose we are given a smooth $1$-parameter family
$\{\Pi_{t}\}_{t\in[0,1]}$ of coupling Poisson tensors on $E$.
For every $t$, $\Pi_{t}$ is a coupling Poisson tensor associated
with geometric data $((\Pi_{t})_{V},\Gamma^{t},F^{t})$, 
where the vertical bivector field $(\Pi_{t})_{V}$,
the Ehresmann connection $\Gamma^{t}$, 
and the $2$-form $F^{t}$ vary smoothly with~$t$. 
In coordinates,
\begin{equation}
\Pi_{t}=-\frac{1}{2}F_{t}^{ij}(\xi,x)\hor_{i}^{t}\wedge\hor_{j}^{t}
+\frac{1}{2}(\Pi_{t})_{V}^{\alpha\beta}(\xi,x)
\frac{\partial}{\partial x^{\alpha}}\wedge 
\frac{\partial}{\partial x^{\beta}}.
\tag{2.46}
\end{equation}
Here we denote $\hor_{i}^{t}=\hor^{\Gamma^{t}}
(\frac{\partial}{\partial\xi^{i}})$ and 
$F_{t}^{is}F_{sj}^{t}=\delta_{j}^{i}$.

Let $Z_{t}$ be a time-dependent vector field on $E$. Then we
have the decomposition
\begin{equation}
Z_{t}=X_{t}+Y_{t},
\tag{2.47}
\end{equation}
where $X_{t}$ and $Y_{t}$ are the horizontal and vertical
components with respect to the connection $\Gamma^{t}$, 
\begin{equation}
X_{t}=X_{t}^{i}(\xi,x)\hor_{i}^{t},\qquad
Y_{t}=Y_{t}^{\sigma}(\xi,x)\frac{\partial}{\partial x^{\sigma}}.
\tag{2.48}
\end{equation}
Let $\partial^{\Gamma^{t}}$ be the covariant exterior
differential (2.2) associated with~$\Gamma^{t}$.

\begin{proposition}
A time-dependent vector field $Z_{t}=X_{t}+Y_{t}$ is a solution
of the equation
\begin{equation}
L_{Z_{t}}\Pi_{t}+\frac{\partial\Pi_{t}}{\partial t}=0
\tag{2.49}
\end{equation}
if and only if the components $X_{t}$ and $Y_{t}$ satisfy the relations
\begin{align}
L_{Y_{t}}\Pi_{V}^{t}+\frac{\partial\Pi_{V}^{t}}{\partial t}&=0,
\tag{2.50}
\\
\partial^{\Gamma^{t}}(\mathbf{i}_{Xt}F^{t})+L_{Y_{t}}F^{t}
+\frac{\partial F^{t}}{\partial t}&=0,
\tag{2.51}
\\
(\Pi^{t})_{V}^{\natural}d(\mathbf{i}_{Xt}F^{t})_{s}
+[\hor_{s}^{t},Y_{t}]
-\frac{\partial\hor_{s}^{t}}{\partial t}&=0.
\tag{2.52}
\end{align}
\end{proposition}

\begin{proof}
First, the straightforward computations give
\begin{align*}
\frac{\partial\Pi_{t}}{\partial t}  
&=\frac{1}{2}F_{t}^{im}
\bigg(\frac{\partial F_{mm^{\prime}}^{t}}{\partial t}\bigg)
F_{t}^{m^{\prime}j}\hor_{i}^{t}\wedge\hor_{j}^{t}
\\
&\qquad
-F_{t}^{ij}\frac{\partial\hor_{i}^{t}}{\partial t}
\wedge\hor_{j}^{t}+\frac{1}{2}\frac{\partial}{\partial t}
(\Pi_{t})_{V}^{\alpha\beta}
\frac{\partial}{\partial x^{\alpha}}\wedge
\frac{\partial}{\partial x^{\beta}}.
\end{align*}
By formulas (2.38), (2.41) and (2.42) we get
\begin{align*}
L_{Z_{t}}\Pi_{t}  
&=\frac{1}{2} F_{t}^{im}[(\partial^{\Gamma^{t}}
(\mathbf{i}_{X_{t}}F^{t}))_{mm^{\prime}}
+L_{Y_{t}}F_{mm^{\prime}}^{t}]
F_{t}^{m^{\prime}j}\hor_{i}^{t}\wedge\hor_{j}^{t}
\\
&\quad 
-F_{t}^{si}\hor_{i}^{t}\wedge(\Pi_{t})_{V}^{\sharp}
(d(X_{t}^{m}F_{ms}^{t}))
-F_{t}^{si}[\hor_{s}^{t},Y_{t}]^{\sigma}\hor_{i}^{t}
\wedge\frac{\partial}{\partial x^{\sigma}}
+L_{Y_{t}}\Pi_{V}.
\end{align*}
Finally, putting these relations into (2.49) and using that
$$
\frac{\partial\hor_{s}^{t}}{\partial t}
\text{ \textit{is a vertical vector field}},
$$
we see that vanishing the terms of degrees $(0,2)$, $(2,0)$, 
and $(1,1)$ leads to Eqs.~(2.50)--(2.52).
\end{proof}

A time-dependent vector field $Z_{t}$ satisfying (2.49) is said
to be an \textit{infinitesimal generator} of the family 
$\{\Pi_{t}\}$ of coupling Poisson tensors on~$E$.

\begin{definition}
We say that two coupling Poisson tensors $\Pi^{\prime}$ and
$\Pi^{\prime\prime}$ on $E$ are \textit{homotopic\/}
if there exists a smooth family $\{\Pi_{t}\}_{t\in[0,1]}$ of
coupling Poisson tensors on $E$  which admits an infinitesimal
generator and joints $\Pi^{\prime}$ and 
$\Pi^{\prime\prime}$, $\Pi_{0}=\Pi^{\prime}$, 
$\Pi_{1}=\Pi^{\prime\prime}$.
\end{definition}

Assume that $Z_{t}$ is an infinitesimal generator of
$\{\Pi_{t}\}$. Let $\Phi_{t}$ be the flow of $Z_{t}$,
\begin{equation}
\frac{d\Phi_{t}}{dt}=Z_{t}\circ\Phi_{t},\qquad 
\Phi_{0}=\id.
\tag{2.53}
\end{equation}
Suppose that for every $t\in[0,1]$ the flow $\Phi_{t}$ is well defined
on an open domain $N_{\flow}\subseteq E$ independent of~$t$.
Then, we have
\begin{equation}
\frac{d}{dt}(\Phi_{t}^{\ast}\Pi_{t})
=\Phi_{t}^{\ast}
\bigg(L_{Z_{t}}\Pi_{t} +\frac{\partial\Pi_{t}}{\partial t}\bigg)=0
\tag{2.54}
\end{equation}
and hence
\begin{equation}
\Phi_{t}^{\ast}\Pi_{t}=\Pi_{0}\qquad\text{on}\quad N_{\flow}.
\tag{2.55}
\end{equation} 
Thus, the family $\{\Pi_{t}\}$ is generated by the flow
$\Phi_{t}$ and the ``initial'' coupling Poisson tensor.

The next question is to formulate some criteria for the
existence of $Z_{t}$ in terms of the geometric data of $\Pi_{t}$
analyzing Eqs.~(2.50)--(2.52). 

\textbf{Sufficient conditions for the existence of $Z_{t}$.}  
Here we assume that we are given a smooth $1$-parameter family
$\{\Pi_{t}\}_{t\in[0,1]}$ of coupling Poisson tensors such that
the vertical part of $\Pi_{t}$ is independent of the parameter~$t$,
\begin{equation}
(\Pi_{t})_{V}=\Upsilon\quad \forall t\in[0,1],
\tag{2.56}
\end{equation}
where $\Upsilon$ is a fiber-tangent Poisson structure. 
Thus, for every $t$, $\Pi_{t}$ is a coupling Poisson tensor
associated with geometric data $(\Upsilon,\Gamma^{t},F^{t})$.

Let us introduce the following notation. 
Denote by $\{,\}_{\Upsilon}$ the Poisson bracket corresponding
to the fiber-tangent Poisson structure $\Upsilon$. 
One can associate to any $Q\in\Omega^{1}(B)\otimes C^{\infty}(E)$
and $\Theta\in\Omega^{k}(B)\otimes C^{\infty}(E)$ an element
$\{Q\wedge\Theta\}_{\Upsilon}\in\Omega^{k+1}(B)\otimes C^{\infty}(E)$ 
defined by
\begin{align}
&\{Q\wedge\Theta\}_{\Upsilon}(u_{0},u_{1},\dots,u_{k})
\tag{2.57}
\\
&\qquad \overset{\text{def}}{=}
\sum_{i=0}^{k}(-1)^{i}
\{Q(u_{i}),\Theta(u_{0},u_{1},\dots,\hat{u}_{i},\dots,u_{k})\}_{\Upsilon}.
\nonumber
\end{align}
In particular,

\begin{itemize}
\item  if $\Theta\in Q\in\Omega^{1}(B)\otimes C^{\infty}(N)$,
\begin{equation}
\{Q\wedge\Theta\}_{\Upsilon}(u_{0},u_{1})
=\{Q(u_{0}),\Theta(u_{1})\}_{\Upsilon}
-\{Q(u_{1}),\Theta(u_{0})\}_{\Upsilon};
\tag{2.58}
\end{equation}

\item  if $\Theta\in\Omega^{2}(B)\otimes C^{\infty}(\Upsilon)$,
\begin{equation}
\{Q\wedge\Theta\}_{\Upsilon}(u_{0},u_{1},u_{2})
=\underset{(u_{0},u_{1},u_{2})}{\mathfrak{S}}
\{Q(u_{0}),\Theta(u_{1},u_{2})\}_{\Upsilon}.
\tag{2.59}
\end{equation}
\end{itemize}

Assume that the data $(\Gamma^{t},F^{t})$ satisfy the following
condition: 
there exists a smooth family of $1$-forms
$Q^{t}\in\Omega^{1}(B)\otimes C^{\infty}(E)$ such that
\begin{align}
\hor^{\Gamma^{t}}(u)
&=\hor^{\Gamma^{0}}(u)+\Upsilon^{\sharp}dQ^{t}(u),
\tag{2.60}
\\
F^{t}&=F^{0}-(\partial^{\Gamma^{0}}Q^{t}
+\frac{1}{2}\{Q^{t}\wedge Q^{t}\}_{\Upsilon})
\tag{2.61}
\end{align}
for every $t\in[0,1]$ and $u\in\mathcal{X}(B)$.

By (2.24), $\hor^{\Gamma^{t}}(u)-\hor^{\Gamma^{0}}(u)$ 
is a vertical Poisson vector field of $\Upsilon$. 
Condition (2.60) says that this difference is a Hamiltonian
vector field with Hamiltonian $Q^{t}(u)$ which must satisfy the
compatibility condition 
\begin{equation}
Q^{0}(u)\in\Casim(E,\Upsilon).
\tag{2.62}
\end{equation}
Furthermore, condition (2.60) implies that the exterior covariant
derivatives $\partial^{\Gamma^{t}}$ and $\partial^{\Gamma^{0}}$
are related by the formula
\begin{equation}
\partial^{\Gamma^{t}}\Theta
= \partial^{\Gamma^{0}}\Theta+\{Q^{t}\wedge \Theta\}_{\Upsilon}.
\tag{2.63}
\end{equation}
In particular,
$$
\partial^{\Gamma^{t}}Q^{t}
=\partial^{\Gamma^{0}}Q^{t}+\{Q^{t}\wedge Q^{t}\}_{\Upsilon}.
$$
It follows that (2.61) can be rewritten in the form
\begin{equation}
F^{t}=F^{0}-\bigg(\partial^{\Gamma^{t}}Q^{t}
-\frac{1}{2}\{Q^{t}\wedge Q^{t}\}_{\Upsilon}\bigg).
\tag{2.64}
\end{equation}

\begin{proposition}
Under assumptions {\rm(2.60), (2.61)}, a time-dependent vector
field $Z_{t}=X_{t}+Y_{t}$ is an infinitesimal generator of the
family $\{\Pi_{t}\}$ if and only if $X_{t}$ and $Y_{t}$ satisfy the
equations 
\begin{align}
L_{Y_{t}}\Upsilon&=0,
\tag{2.65}
\\
\partial^{\Gamma^{t}}
\bigg(\frac{\partial Q^{t}}{\partial t}-\mathbf{i}_{Xt} F^{t}\bigg)
-L_{Y_{t}}F^{t}&=0,
\tag{2.66}
\\
\Upsilon^{\sharp}d\bigg(\frac{\partial Q^{t}(u)}{\partial t}
-\mathbf{i}_{Xt} F^{t}(u)\bigg)
+[Y_{t},\hor^{t}(u)]&=0.
\tag{2.67}
\end{align}
\end{proposition}

\begin{proof}
By (2.61) and (2.63) we deduce
$$
\frac{\partial F^{t}}{\partial t}
=-\bigg[\partial^{\Gamma^{0}}
\bigg(\frac{\partial Q^{t}}{\partial t}\bigg)
+\{Q^{t}\wedge\frac{\partial Q^{t}}{\partial t}\}_{\Upsilon}\bigg]
=-\partial^{\Gamma^{t}}
\bigg(\frac{\partial Q^{t}}{\partial t}\bigg).
$$

Condition (2.60) implies
\begin{equation}
\frac{\partial\hor_{i}^{^{t}}}{\partial t}=\Upsilon^{\sharp}dQ^{t}.
\tag{2.68}
\end{equation}
Putting these relations into (2.51) and (2.52) gives (2.66) and (2.67).
\end{proof}

\begin{corollary}
We have the variation of parameter formula for~$\Pi_{t}$
\begin{align}
\frac{\partial\Pi_{t}}{\partial t}  
& =-\frac{1}{2}F_{t}^{im}\bigg(\partial^{\Gamma^{t}}
\bigg(\frac{\partial Q^{t}}{\partial t}\bigg)\bigg)_{mm^{\prime}}
F_{t}^{m^{\prime}j}
\hor_{i}^{t}\wedge\hor_{j}^{t}
\tag{2.69}
\\
&\qquad 
+F_{t}^{ij}\hor_{j}^{t}\wedge\Upsilon^{\natural}
\bigg(d\bigg(\frac{\partial Q_{i}^{t}}{\partial t}\bigg)\bigg).
\nonumber
\end{align}
\end{corollary}

Now, choosing $Y_{t}=0$ in (2.66), (2.67), we arrive at the key
observation. 

\begin{theorem}
Let $\{\Pi_{t}\}$ be a family of coupling Poisson structures
associated with geometric data $(\Upsilon,\Gamma^{t},F^{t})$. If
conditions {\rm(2.60), (2.61)} hold, then $\{\Pi_{t}\}$ admits
the infinitesimal generator $X_{t}=X_{t}^{i}(\xi,x)\hor_{i}^{t}$
which is uniquely determined by 
\begin{equation}
\mathbf{i}_{Xt}F^{t}=\frac{\partial Q^{t}}{\partial t}.
\tag{2.70}
\end{equation}
\end{theorem}

Remark that $X_{t}$ is well defined by (2.70) because of the
nondegeneracy property (2.20) for $F^{t}$.

Let $Y_{t}=\Upsilon^{\natural}dh_{t}$ be the vertical Hamiltonian
time-dependent vector relative to $\Upsilon$ with a
time-dependent function $h_{t}\in C^{\infty}(N)$. 
By (2.24), the horizontal lift $\hor^{\Gamma^{t}}(u)$ is a
vertical Poisson vector field of $\Upsilon$ and hence 
\begin{equation}
[\hor^{\Gamma^{t}}(u),\Upsilon^{\natural}dh_{t}]
=\Upsilon^{\natural}d(L_{\hor^{\Gamma^{t}}(u)}h_{t}).
\tag{2.71}
\end{equation}
It follows from (2.65)--(2.67), that under conditions (2.60),
(2.61), finding solutions of (2.49) in the form
$Z_{t}=X_{t}+\Upsilon^{\natural}dh_{t}$ is 
reduced to solving the following two equations:
\begin{align}
&\partial^{\Gamma^{t}}
\bigg(\frac{\partial Q^{t}}{\partial t}-\mathbf{i}_{X_{t}}F^{t}\bigg)
=\{h_{t},F^{t}\}_{\Upsilon},
\tag{2.72}
\\
&\Upsilon^{\sharp}d\bigg(\frac{\partial Q^{t}(u)}{\partial t}
-\mathbf{i}_{X_{t}}F^{t}(u)-L_{\hor^{\Gamma^{t}}(u)}h_{t}\bigg)
=0.
\tag{2.73}
\end{align}

\textbf{Necessary conditions.} 
Theorem 2.18 says that conditions
(2.60), (2.61) are sufficient for the existence of an
infinitesimal generator of $\{\Pi_{t}\}$. 
Now we assume that only (2.60) holds. In general, the
implication (2.60) $\Longrightarrow$ (2.61) is not true. By
(2.26) one can conclude that (2.61) holds only modulo an element
in $\Omega^{2}(B)\otimes\Casim(E,\Upsilon)$.
The next statement shows that condition (2.61) is necessary for
the existence on an infinitesimal generator 
of $\{\Pi_{t}\}$ in the class of vector fields of the form 
$Z_{t}=X_{t}+\Upsilon^{\natural}dh_{t}$.

\begin{theorem}
Let $\{\Pi_{t}\}$ be a family of coupling Poisson structures
associated with geometric data $(\Upsilon,\Gamma^{t},F^{t})$.
Assume that 

{\rm(a)} condition {\rm(2.60)} holds for a certain family
$\{Q^{t}\}${\rm;}

{\rm(b)} $\{\Pi_{t}\}$ admits an infinitesimal generator of the
form 
\begin{equation}
Z_{t}=X_{t}+\Upsilon^{\natural}dh_{t},
\tag{2.74}
\end{equation}
for $h^{t}\in C^{\infty}(E)$.

Then, the family $\{Q^{t}\}$ can be chosen such that {\rm(2.61)}
holds. 
\end{theorem}

\begin{proof}
Suppose $\{Q^{t}\}$ in (2.60) and an infinitesimal generator
$Z_{t}$ in (2.74) are fixed. Remark that $\{Q^{t}\}$ is
unique up to the transformation $Q^{t}\mapsto Q^{t}+K^{t}$, 
where $K^{t}\in\Omega^{1}(B)\otimes\Casim(E,\Upsilon)$. 
Our point is to choose $K^{t}$ to satisfy (2.61). By 
Proposition~2.14, $Z_{t}$ must satisfy Eqs.~(2.51) and~(2.52). 
First, let us consider (2.52)  for $Y_{t}=\Upsilon^{\natural}dh_{t}$.
Using the identity
$$
(\partial^{\Gamma^{t}}h^{t})(u)=L_{\hor^{t}(u)}h^{t}\quad
\forall u\in\mathcal{X}(B),
$$
by (2.71) we deduce
$$
[\hor_{s}^{t},Y_{t}]=\Upsilon^{\natural}d(\partial^{\Gamma^{t}}h^{t})_{s},
$$
where $\partial^{\Gamma^{t}}h^{t}\in\Omega^{1}(B)\otimes C^{\infty}(E)$.
Putting this and (2.68) into (2.52) gives
$$
\Upsilon^{\natural}d\bigg(  
(\mathbf{i}_{Xt}F^{t})_{s}+(\partial^{\Gamma^{t}}h^{t})_{s}
-\frac{\partial Q_{s}^{t}}{\partial t}\bigg)=0.
$$
Consequently
$$
(\mathbf{i}_{Xt}F^{t})=\frac{\partial Q^{t}}{\partial t}
-\partial^{\Gamma^{t}}h^{t}+\beta^{t}
$$
for a certain $\beta^{t}\in\Omega^{1}(B)\otimes\Casim(E,\Upsilon)$. 
Applying the operator $\partial^{\Gamma^{t}}$ to both sides of
this equality, we have
\begin{align}
\partial^{\Gamma^{t}}(\mathbf{i}_{Xt}F^{t})  
& =\partial^{\Gamma^{t}}
\bigg(\frac{\partial Q^{t}}{\partial t}
-\partial^{\Gamma^{t}}h^{t}+\beta^{t}\bigg)
\tag{2.75}\\
& =\partial^{\Gamma^{t}}
\bigg(\frac{\partial Q^{t}}{\partial t}\bigg)
-(\partial^{\Gamma^{t}})^{2}h^{t}+\partial^{\Gamma^{t}}\beta^{t}
\nonumber
\end{align}
Next, consider Eq.~(2.51).  
Remark that we have the identity (see (2.79):
\begin{equation}
L_{Y_{t}}F^{t}=\{h^{t},F^{t}\}_{\Upsilon}
=(\partial^{\Gamma^{t}})^{2}h^{t}.
\tag{2.76}
\end{equation}
By using (2.75) and (2.76), we rewrite (2.51) in the form
\begin{equation}
\partial^{\Gamma^{t}}
\bigg(\frac{\partial Q^{t}}{\partial t}\bigg)
+\partial^{\Gamma^{T}}\beta^{t}+\frac{\partial F^{t}}{\partial t}=0.
\tag{2.77}
\end{equation}
Now taking into account (2.63) and that 
$\partial^{\Gamma^{t}}\beta^{t}=\partial^{\Gamma^{0}}\beta^{t}$, 
from (2.77) we deduce
\begin{align*}
\frac{\partial F^{t}}{\partial t}  
& =-\bigg(\partial^{\Gamma^{0}}
\bigg(\frac{\partial Q^{t}}{\partial t}\bigg)
+\bigg\{Q^{t},\frac{\partial Q^{t}}{\partial t}\bigg\}_{\Upsilon}
\bigg)
-\partial^{\Gamma^{0}}\beta^{t}\\
& =-\frac{\partial}{\partial t}
\bigg[\partial^{\Gamma^{0}}Q^{t}
+\frac{1}{2}\{Q^{t}\wedge Q^{t}\}_{\Upsilon}\bigg]
-\partial^{\Gamma^{0}}\beta^{t}.
\end{align*}
Finally, integrating of this identity gives
\begin{equation}
F^{t}=F_{0}-\bigg[\partial^{\Gamma^{0}}Q^{t}+\frac{1}{2}\{Q^{t}\wedge
Q^{t}\}_{\Upsilon}\bigg]
-\partial^{\Gamma^{0}}\bigg(\int_{0}^{t}\beta^{\tau}\,d\tau\bigg).
\tag{2.78}
\end{equation}
Thus, after the transformation
$$
Q^{t}\mapsto Q^{t}+\int_{0}^{t}\beta^{t}dt,
$$
(2.61) is satisfied.
\end{proof}

\subsection{Relative Casimir $2$-cocycles}
Denote by $\Coup(E,\Upsilon)$ the set of all coupling Poisson
tensors $\Pi$ with fixed vertical part, $\Pi_{V}=\Upsilon$. Here
we introduce an invariant of a pair of bivector fields in
$\Coup(E,\Upsilon)$ which takes values in the space 
$\Omega^{2}(B)\otimes\Casim (E,\Upsilon)$. This invariant is
well defined for coupling Poisson tensors
whose geometric data belong to the same equivalence class.

Pick $\Pi\in\Coup_{B}(E,\Upsilon)$ associated with geometric
data $(\Upsilon,\Gamma,F)$. Consider the covariant exterior
derivative $\partial^{\Gamma}$. By (2.3) and (2.26) we have the
identity 
\begin{align}
& ((\partial^{\Gamma})^{2}\Theta)
(u_{0},u_{1},\dots,u_{i},\dots,u_{j},\dots,u_{k+1})
\tag{2.79}
\\
&\qquad 
=\sum_{0\leq i<j\leq k+1}(-1)^{i+j}
\{F(u_{i},u_{j}),\Theta(u_{0},u_{1},\dots,\hat{u}_{i},\dots,
\hat{u}_{j},\dots,u_{k+1})\}_{\Upsilon}.
\nonumber
\end{align}

We also need the following formulas which can be easily derived
from definitions. For any $Q\in\Omega^{1}(B)\otimes
C^{\infty}(E)$ we have 
\begin{align}
(\partial^{\Gamma})^{2}Q&=\{Q\wedge F\}_{\Upsilon},
\tag{2.80}
\\
\partial^{\Gamma}\{Q\wedge Q\}_{\Upsilon}
&=-2\{Q\wedge\partial^{\Gamma}Q\}_{\Upsilon},
\tag{2.81}
\\
\{Q\wedge\{Q\wedge Q\}_{\Upsilon})\}_{\Upsilon}&=0.
\tag{2.82}
\end{align}

Note that the Lie derivative $L_{\hor_{u}^{\Gamma}}$ preserves
the subspace $\Casim(E,\Upsilon)\subset C^{\infty}(E)$. This
means that the restriction of $\partial^{\Gamma}$ to 
$\Omega^{k}(B)\otimes\Casim(E,\Upsilon)$ is well defined. 
We will denote this restriction by
\begin{equation}
\partial_{0}^{\Gamma}:\,\Omega^{k}(B)\otimes\Casim(E,\Upsilon)
\rightarrow\Omega^{k+1}(B)\otimes\Casim(E,\Upsilon).
\tag{2.83}
\end{equation}

It follows from (2.79) that $\partial_{0}^{\Gamma}$ is a
\textit{coboundary operator},
$$
(\partial_{0}^{\Gamma})^{2}=0.
$$

Now let us take another coupling Poisson tensor 
$\tilde{\Pi}\in \Coup(E,\Upsilon)$ associated with geometric
data $(\Upsilon,\tilde{\Gamma},\tilde{F})$. 
Assume that the pair $(\tilde{\Pi},\Pi)$
satisfies the condition: there exists $Q\in\Omega^{1}(B)\otimes
C^{\infty}(E)$ such that
\begin{equation}
\hor^{\tilde{\Gamma}}(u)=\hor^{\Gamma}(u)+\Upsilon^{\sharp}dQ(u)
\tag{2.84}
\end{equation}
for $u\in\mathcal{X}(B)$. Then, for every
$\Theta\in\Omega^{k}(B)\otimes C^{\infty}(E)$
\begin{equation} 
\partial^{\tilde{\Gamma}}\Theta
=\partial^{\Gamma}\Theta+\{Q\wedge \Theta\}_{\Upsilon}.
\tag{2.85}
\end{equation}
This implies
\begin{equation}
\partial_{0}^{\Gamma}=\partial_{0}^{\tilde{\Gamma}}.
\tag{2.86}
\end{equation}
Let us associate to the pair $(\tilde{\Pi},\Pi)$ the following $2$-form
\begin{equation}
C=C_{\tilde{\Pi}\Pi}\overset{\text{def}}{=}
\tilde{F}-F
+\bigg(\partial^{\Gamma}Q+\frac{1}{2}\{Q\wedge Q\}_{\Upsilon}\bigg).
\tag{2.87}
\end{equation}
Then,
\begin{equation}
C\in\Omega^{2}(B)\otimes\Casim(E,\Upsilon).
\tag{2.88}
\end{equation}
Indeed, by (2.81) we get
$$
\Curv^{\tilde{\Gamma}}=\Curv^{\Gamma}
-\Upsilon^{\sharp}d(\partial^{\Gamma}Q
+\frac{1}{2}\{Q\wedge Q\}_{\Upsilon}).
$$
Comparing this with
$$
\Curv^{\tilde{\Gamma}}=\Upsilon^{\sharp}d\tilde{F}
\qquad\text{and}\qquad
\Curv^{\Gamma}=\Upsilon^{\sharp}dF
$$
leads to (2.85).

\begin{proposition}
The $2$-form $C$ in {\rm(2.87)} is a $2$-cocycle,
\begin{equation}
\partial_{0}^{\Gamma}C=0.
\tag{2.89}
\end{equation}
\end{proposition}

\begin{proof}
It follows from (2.25), (2.84) that
$$
\partial^{\Gamma}\tilde{F}
=(\partial^{\Gamma}-\partial^{\tilde{\Gamma}})\tilde{F}
=-\{Q\wedge\tilde{F}\}_{\Upsilon}.
$$
On the other hand, by using (2.82) and (2.88) we derive
\begin{align*}
\{Q\wedge\tilde{F}\}_{\Upsilon}  
& =\{Q\wedge F\}_{\Upsilon}+\{Q\wedge C\}_{\Upsilon}
-\{Q\wedge\partial^{\Gamma}Q\}_{\Upsilon}
-\frac{1}{2}\{Q\wedge\{Q\wedge Q\}_{\Upsilon})\}_{\Upsilon}
\\
& =\{Q\wedge F\}_{\Upsilon}-\{Q\wedge\partial^{\Gamma}Q\}_{\Upsilon}.
\end{align*}
Thus,
$$
\partial^{\Gamma}\tilde{F}
=-\{Q\wedge\partial^{\Gamma}Q\}_{\Upsilon}+\{Q\wedge F\}_{\Upsilon}.
$$
Finally, taking into account (2.80),(2.81), we get
\begin{align*}
\partial^{\Gamma}C  
& =\partial^{\Gamma}\tilde{F}+(\partial^{\Gamma})^{2}Q
+\frac{1}{2}\partial^{\Gamma}\{Q\wedge Q\}_{\Upsilon}
\\
& =-\{Q\wedge F\}_{\Upsilon}+\{Q\wedge\partial^{\Gamma}Q\}_{\Upsilon}
+(\partial^{\Gamma})^{2}Q+\frac{1}{2}\partial^{\Gamma}
\{Q\wedge Q\}_{\Upsilon}
\\
& =0.
\end{align*}
\end{proof}

If (2.84) holds, then we shall write $\tilde{\Pi}\sim_{Q}\Pi$.
It is clear that condition (2.84) defines an \textit{equivalence
relation} for the elements of $\Coup(E,\Upsilon)$.

\begin{definition}
For any $\tilde{\Pi},\Pi\in\Coup(E,\Upsilon)$ such that
$\tilde{\Pi}\sim_{Q}\Pi$, the $2$-form $C=C_{\tilde{\Pi}\Pi}$ in
(2.87) will be called the \textit{Casimir $2$-cocycle} of the
pair $(\tilde{\Pi},\Pi)$. 
\end{definition}

Recall that $Q$ is uniquely determined by (2.84) up to the
transformation 
$$
Q\mapsto Q+K,
$$
for arbitrary $K\in\Omega^{1}(B)\otimes\Casim(E,\Upsilon)$.
Then $C$ is transformed by the rule
$$
C\mapsto C+\partial_{0}^{\Gamma}K,
$$
which preserves the $\partial_{0}^{\Gamma}$-\textit{cohomology
class} $[C]=[C_{\tilde{\Pi}\Pi}]$ of $C=C_{\tilde{\Pi}\Pi}$. 
Thus, $[C_{\tilde{\Pi}\Pi}]$ is an invariant of the pair
$(\tilde{\Pi},\Pi)$. Here are some properties. 

\begin{proposition}
For any coupling Poisson tensors $\Pi$, $\Pi^{\prime}$,
$\Pi^{\prime\prime}$ from the same equivalence class, 
the following identities hold:
\begin{align}
[C_{\Pi^{\prime}\Pi}]&=-[C_{\Pi\Pi^{\prime}}],
\tag{2.90}
\\
[ C_{\Pi^{\prime}\Pi}]+[C_{\Pi\Pi^{\prime\prime}}]
+[C_{\Pi^{\prime \prime}\Pi^{\prime}}]&=0.
\tag{2.91}
\end{align}
\end{proposition}

\begin{proof}
Below we perform computations modulo the elements in 
$\Omega^{1}(B)\otimes\Casim(E,\Upsilon)$. By (2.85), we have
\begin{align*}
C_{\Pi^{\prime}\Pi}  
& =F-F^{\prime}
-\bigg(\partial^{\Gamma^{\prime}}Q-\frac{1}
{2}\{Q\wedge Q\}_{\Upsilon}\bigg)
\\
& =F-F^{\prime}
-\bigg(\partial^{\Gamma}Q+\frac{1}{2}\{Q\wedge Q\}_{\Upsilon}\bigg)
\\
& =-C_{\Pi\Pi^{\prime}}.
\end{align*}
Moreover, (2.84) implies
$$
Q^{\Gamma^{\prime\prime}\Gamma}
=Q^{\Gamma^{\prime\prime}\Gamma^{\prime}}+Q^{\Gamma^{\prime}\Gamma}
$$
and then
\begin{align*}
C_{\Pi^{\prime\prime}\Pi}
&\overset{\text{def}}{=}
F^{\prime\prime}-F+(\partial^{\Gamma} Q^{\Gamma^{\prime\prime}\Gamma}
+\frac{1}{2}\{Q^{\Gamma^{\prime\prime}\Gamma}\wedge 
Q^{\Gamma^{\prime\prime}\Gamma}\}_{\Upsilon}
\\
&=(F^{\prime\prime}-F^{\prime})+(F^{\prime}-F)
+\partial^{\Gamma}
Q^{\Gamma^{\prime\prime}\Gamma^{\prime}}
+\partial^{\Gamma}Q^{\Gamma^{\prime}\Gamma}
\\
&\qquad
+\frac{1}{2}\{Q^{\Gamma^{\prime\prime}\Gamma^{\prime}}
\wedge 
Q^{\Gamma^{\prime\prime}\Gamma^{\prime}} \}_{\Upsilon}
+\frac{1}{2}\{Q^{\Gamma^{\prime}\Gamma}\wedge 
Q^{\Gamma^{\prime}\Gamma} \}_{\Upsilon}
+\{Q^{\Gamma^{\prime\prime}\Gamma^{\prime}}
\wedge Q^{\Gamma^{\prime}\Gamma} \}_{\Upsilon}.
\end{align*}
Finally, using again (2.85),
$$
\partial^{\Gamma}Q^{\Gamma^{\prime\prime}\Gamma^{\prime}}
=\partial^{\Gamma^{\prime}}Q^{\Gamma^{\prime\prime}\Gamma^{\prime}}
-\{Q^{\Gamma^{\prime}\Gamma}\wedge 
Q^{\Gamma^{\prime\prime}\Gamma}\}_{\Upsilon},
$$
we derive (2.91).
\end{proof}

Let $\Ham(E,\Upsilon)$ be the Lie algebra of Hamiltonian vector
fields of the fiber-tangent Poisson tensor $\Upsilon$. Denote
by $\Poiss_{V}(E,\Upsilon)$ the Lie algebra of all vertical
Poisson vector fields of $\Upsilon$. 
Then $\Ham(E,\Upsilon)$ is an ideal of $\Poiss_{V}(E,\Upsilon)$. 
Introduce the first reduced Poisson cohomology space as
\begin{equation} 
H_{V}^{1}(E,\Upsilon)
=\frac{\Poiss_{V}(E,\Upsilon)}{\Ham(E,\Upsilon)}.
\tag{2.92}
\end{equation}
The following statement is evident.

\begin{proposition}
Assume that
\begin{equation}
H_{V}^{1}(E,\Upsilon)=\{0\}.
\tag{2.93}
\end{equation}
Then for arbitrary $\tilde{\Pi},\Pi\in\Coup(E,\Upsilon)$ there
exists a $Q\in\Omega^{1}(B)\otimes C^{\infty}(E)$ such that
$\tilde{\Pi} \sim_{Q}\Pi$, and hence the Casimir cocycle 
$[C_{\tilde{\Pi}\Pi}]$ is well defined.
\end{proposition}

Under condition (2.93), one can associated to $\Upsilon$ an
\textit{intrinsic coboundary operator}
\begin{equation}
\partial_{0}:\,\Omega^{k}(B)\otimes\Casim(E,\Upsilon)
\rightarrow\Omega^{k+1}(B)\otimes\Casim(E,\Upsilon),
\tag{2.94}
\end{equation}
given by $\partial_{0}=\partial_{0}^{\Gamma}$, where $\Gamma$ is
a Poisson connection on $(E,\Gamma)$. It is clear that this
definition is independent of $\Gamma$. 
The properties of $\partial_{0}$ are controlled by the
(singular) symplectic foliation of~$\Upsilon$.

\textbf{Inner Poisson automorphisms.} 
Denote by $\Aut(E,\Upsilon)$ the group of all fiber preserving
Poisson automorphisms of $(E,\Upsilon)$. 
By $\Inn(E,\Upsilon)$ we denote the subgroup of
$\Aut(E,\Upsilon)$ corresponding to \textit{inner Poisson 
automorphisms}. Recall that a diffeomorphism $g:\,E\rightarrow E$
is called an \textit{inner Poisson automorphism}
\cite{GiGo,Fe1}
if there exists a smooth family of Hamiltonian functions 
$h_{t}:\,E\rightarrow\mathbb{R}$ such that the flow
$g_{t}$ of the time-dependent Hamiltonian vector field
$\mathcal{V}^{\sharp}d(h^{t})$ joints $g$ with the identity map,
$g_{0}=\id$ and $g_{1}=g$.

By Proposition 2.9, the groups $\Aut(E,\Upsilon)$ and
$\Inn(E,\Upsilon)$ act naturally on $\Coup(E,\Upsilon)$.

\begin{proposition}
Let $\tilde{\Pi},\Pi\in\Coup(E,\Upsilon)$. 
If $\tilde{\Pi}\sim_{Q}\Pi$, then for every 
$g\in\Inn(E,\Upsilon)$, the push-forward 
$g_{\ast}\Pi\in\Coup(E,\Upsilon)$ is a coupling
Poisson tensor associated with geometric data 
$(\Upsilon,g_{\ast}\Gamma,g_{\ast}F)$ and such that
\begin{align}
\tilde{\Pi}&\sim_{Q^{1}}g_{\ast}\Pi,
\tag{2.95}
\\
[ C_{\Pi,\tilde{\Pi}}]&=[C_{g_{\ast}\Pi,\tilde{\Pi}}]
\tag{2.96}
\end{align}
for a certain $Q^{1}\in\Omega^{1}(B)\otimes C^{\infty}(E)$.
\end{proposition}

\begin{proof} 
We use the flow $g_{t}$ of the time-dependent Hamiltonian vector
field $\mathcal{V}^{\sharp}d(h^{t})$ to define a family of
coupling Poisson tensors $\{\Pi_{t}\}_{t\in[0,1]}$ by
\begin{equation}
\Pi_{t}=(g_{t})_{\ast}\Pi.
\tag{2.97}
\end{equation}
By Proposition 2.9, the corresponding geometric data are
$$
(\Upsilon,\Gamma^{t},F^{t})
=\big(\Upsilon,(g_{t})_{\ast}\Gamma,(g_{t})_{\ast}F\big),
$$
where $(\Upsilon,\Gamma,F)$ are geometric data of $\Pi$. Then, 
$\Pi_{1}=\tilde{\Pi}$ and $\Pi_{0}=\Pi$. 
Taking into account (2.71) for $\hor^{\Gamma^{t}}(u)$, we get
\begin{align*}
\hor^{\Gamma^{t}}(u)-\hor^{\Gamma}(u)  
& =\int_{0}^{t}[\hor^{\tau}(u),\mathcal{V}^{\sharp}d(h_{\tau})]\,d\tau
\\
& =\mathcal{V}^{\sharp}
d\bigg(\int_{0}^{t}L_{\hor^{\tau}(u)}h_{\tau}\,d\tau\bigg).
\end{align*}
Thus
$$
\hor^{\tilde{\Gamma}}(u)=\hor^{\Gamma^{t}}(u)+\mathcal{V}^{\sharp}dQ^{t},
$$
where
$$
Q^{t}=Q+\int_{0}^{t}L_{\hor^{\tau}(u)}h_{\tau}\,d\tau.
$$
In particular, at $t=1$ this gives (2.95).

Now applying Theorem 2.19 to family (2.97) and 
$Z_{t}=\mathcal{V}^{\sharp}d(h^{t})$ shows that 
$\beta^{\tau}=0$ in (2.78) and hence $C_{\Pi_{t},\Pi}=0$.
From here and property (2.91) we deduce (2.96).
\end{proof}

\begin{corollary}
The orbit of the action of $\Inn(E,\Upsilon)$ through every
$\Pi\in\Coup(E,\Upsilon)$ belongs to the equivalence class of
$\Pi$. Moreover,
\begin{equation}
C_{g_{\ast}\Pi,\Pi}=0
\tag{2.98}
\end{equation}
for every $g\in\Inn(E,\Upsilon)$.
\end{corollary}

\section{Poisson equivalence}

Here, using the results of previous sections, we obtain some
criteria for the equivalence of Poisson structures in the
neighborhood of an entire common symplectic leaf.

Let $\pi:\,E\rightarrow B$ be a vector bundle over a connected
base $B$ equipped with symplectic structure $\omega$. We shall
identify $B$ with the image of the zero section.

\begin{proposition}
Let $\Pi$ be a Poisson tensor on $E$ with property:
\begin{equation}
(B,\omega)\text{ \textit{is a symplectic leaf of} }\Pi.
\tag{3.1}
\end{equation}
Then, there exists a neighborhood $N$ of $B$ in $E$ such that
$\Pi$ is a coupling Poisson tensor on $N$ associated to
geometric data $(\Pi_{V},\Gamma,F)$ satisfying the following
conditions: 
\begin{align}
\rank\Pi_{V}&=0\qquad\text{on}\quad B,
\tag{3.2}
\\
\hor^{\Gamma}(u) |_{B}&=u,
\tag{3.3}
\\
\pi^{\ast}F |_{B}&=\omega
\tag{3.4}
\end{align}
for every $u\in\mathcal{X}(B)$.
\end{proposition}

\begin{proof}
We have to show that $\Pi$ is horizontally nondegenerate in a
neighborhood of $B$. Let $\mathbb{V}\subset TE$ be the vertical
subbundle and $\mathbb{V}^{0}\subset T^{\ast}E$ be the
annihilator of $\mathbb{V}$. Consider the 
plane distribution on $E$ associated with $\Pi$:
\begin{equation}
\mathbb{H}_{m}=\Pi_{m}^{\sharp}(\mathbb{V}_{m}^{0}),\qquad
m\in E.
\tag{3.5}
\end{equation}
For every $\xi\in B$ we have the decomposition
$$
T_{\xi}E=T_{\xi}B\oplus E_{\xi},
$$
where $\mathbb{V}_{\xi}=E_{\xi}$. Then,
$$
T_{\xi}^{\ast}E=(T_{\xi}B)^{0}\oplus\mathbb{V}_{\xi}^{0}.
$$
On the other hand, it follows from (3.1) that
$$
\Pi_{\xi}^{\sharp}(T_{\xi}^{\ast}E)=T_{\xi}B,
$$
which says that $\Pi_{\xi}^{\sharp}((T_{\xi}B)^{0})=\{0\}$. 
Consequently, 
\begin{equation}
\mathbb{H}_{\xi}=T_{\xi}B
\tag{3.6}
\end{equation}
for every $\xi\in B$. Thus, $\mathbb{H}$ is a complementary of
$\mathbb{V}$ at every point of $B$. This means that there exists
an open neighborhood $N$ of $B$ in $E$ such that
\begin{equation}
T_{N}E=\mathbb{H}\oplus\mathbb{V}.
\tag{3.7}
\end{equation}
Properties (3.3) and (3.4) follow from (3.1) and (3.6). Remark that
$\rank\Pi=\rank\Pi_{H}+\rank\Pi_{V}$. But 
$\rank\Pi=\rank\Pi_{H}=\dim B$ at $\xi\in B$. This proves (3.2).
\end{proof}

Denote by $\mathfrak{Coup}_{B}(E)$ the set of germs at $B$ of all
Poisson structures $\Pi$ on $E$ satisfying condition (3.1). Let
$\phi$ be a diffeomorphism between two neighborhoods of $B$
which is identical on $B$, $\phi |_{B}=\id_{B}$. Then, for
every $\Pi\in\mathfrak{Coup}_{B}(E)$ we have
$\phi_{\ast}\Pi\in\mathfrak{Coup}_{B}(E)$.

We are interested in the equivalence between the coupling Poisson
tensors in $\mathfrak{Coup}_{B}(E)$ in the following class of
diffeomorphisms identical on~$B$.

Denote by $\mathfrak{Diff}_{B}^{0}(E)$ the group of germs at $B$ of
diffeomorphisms $\phi$ satisfying the condition: 
in a neighborhood of $B$ in $E$ the exists a time-dependent
vector field $Z_{t}$ ($t\in[0,1]$) such that
\begin{equation}
Z_{t} |_{B}=0
\tag{3.8}
\end{equation}
and the flow $\Phi_{t}$ of $Z_{t}$ joints $\phi$ and the
identity map, $\Phi_{1}=\phi$,  $\Phi_{0}=\id$.

Introduce also the subgroup $\mathfrak{End}_{B}^{0}(E)$ in 
$\mathfrak{Diff}_{B}^{0}(E)$ consisting of germs of diffeomorphisms
$g$ satisfying the condition: then 
there exists a time-dependent \textit{vertical} vector field
$Y_{t}$ such that
\begin{equation}
Y_{t} |_{B}=0
\tag{3.9}
\end{equation}
and the flow of $Y_{t}$ joints $g$ and the identity map. In
particular, every $g\in\mathfrak{End}_{B}^{0}(E)$ is a fiber
preserving diffeomorphism identical on $B$, 
$\pi\circ g=\pi$ and $g |_{B}=\id_{B}$.

\begin{theorem}
Let $\Pi,\tilde{\Pi}\in\mathfrak{Coup}_{B}(E)$ be two coupling
Poisson tensors associated with geometric data
$(\Pi_{V},\Gamma,F)$ and $(\tilde{\Pi}_{V},\tilde{\Gamma},\tilde{F})$, 
respectively. If $\Pi$ and $\tilde{\Pi}$ are
isomorphic by a diffeomorphism $\phi\in\mathfrak{Diff}_{B}^{0}(E)$,
\begin{equation}
\phi^{\ast}\tilde{\Pi}=\Pi,
\tag{3.10}
\end{equation}
then the vertical Poisson structures $\Pi_{V}$ and
$\tilde{\Pi}_{V}$ are isomorphic by a $g\in\mathfrak{End}_{B}^{0}(E)$,
\begin{equation}
g^{\ast}\tilde{\Pi}_{V}=\Pi_{V}.
\tag{3.11}
\end{equation}
\end{theorem}

\begin{proof}
Let $Z_{t}$ be the time-dependent vector field generating the
diffeomorphism $\phi$. By (3.8) there exists a neighborhood
$N_{\flow}$ of $B$ in $E$ independent of $t$ such that the flow
$\Phi_{t}$ of $Z_{t}$ is well defined on $N_{\flow}$ for all 
$t\in[0,1]$, $\Phi_{t}:\,N_{\flow}\rightarrow E$. It is clear that
$$
\Phi_{t} |_{B}=\id_{B},
$$
$\Phi_{1}=\phi$ and $\Phi_{0}=\id$. For every $t\in[0,1]$,
one can define the Poisson tensor
$$
\Pi_{t}=(\Phi_{t})_{\ast}\Pi
$$
satisfying (3.1). Then, by Proposition 3.1, $\Pi_{t}$ is a
coupling Poisson tensor associated with geometric data
$(\Pi_{V}^{t},\Gamma^{t},F^{t})$ varying smoothly with $t$. 
In particular,
$$
(\Pi_{V}^{0},\Gamma^{0},F^{0})=(\Pi_{V},\Gamma,F)
$$
and
$$
(\Pi_{V}^{1},\Gamma^{1},F^{1})=(\tilde{\Pi}_{V},\tilde{\Gamma},\tilde{F}).
$$
Thus, $Z_{t}$ is an infinitesimal generator of the smooth family
of coupling Poisson structures $\{\Pi_{t}\}_{t\in[0,1]}$, that is, 
$Z_{t}$ satisfies (2.50). For every $t$, we have the decomposition
\begin{equation}
Z_{t}=X_{t}+Y_{t}
\tag{3.12}
\end{equation}
into the horizontal $X_{t}$ and vertical $Y_{t}$ parts with
respect to the connection $\Gamma^{t}$. It follows from (2.50)
and Proposition~2.14 that the time-dependent vertical vector 
field $Y_{t}$ must satisfy the equation
$$
L_{Y_{t}}\Pi_{V}^{t}+\frac{\partial\Pi_{V}^{t}}{\partial t}=0.
$$
Furthermore, by (3.8) and (3.12) we have
$$
Y_{t} |_{B}=0\quad\forall t.
$$
Shrinking the neighborhood $N_{\flow}$, we can arrange that the
flow $g_{t}$ of $Y_{t}$ is well defined on $N_{\flow}$. Then,
$g_{t}^{\ast}\Upsilon_{t}=\Upsilon$ and hence
$$
g^{\ast}\tilde{\Upsilon}=\Upsilon,\qquad g=g_{t} |_{t=1}.
$$
\end{proof}

If (3.10) holds for a certain $\phi\in\mathfrak{Diff}_{B}^{0}(E)$, 
then we say that the coupling Poisson structures 
$\Pi$ and $\tilde{\Pi}$ are \textit{equivalent}. 
Moreover, we can introduce the equivalence relation
for fiber-tangent Poisson structures on $E$ saying that two
fiber-tangent Poisson structures vanishing on $B$ are equivalent
if they are isomorphic by a diffeomorphism in $\mathfrak{End}_{B}^{0}(E)$.
Then,
Theorem~3.2 can be reformulated in the following manner:
\textit{the vertical components of equivalent coupling Poisson
tensors in} $\mathfrak{Coup}_{B}(E)$ \textit{are also equivalent}.

Combining the Local Splitting Theorem and Theorem~3.2, we arrive
at the fact that the vertical part of every coupling Poisson
tensor in $\mathfrak{Coup}_{B}(E)$ is locally trivial (in the sense
of Definition~2.8). 

\begin{proposition}
Let $\Pi\in\mathfrak{Coup}_{B}(E)$ be a coupling Poisson structure,
and let $\Pi_{V}$ be the vertical part. For every $\xi\in B$
there exist a neighborhood $U$ of $\xi$ in $B$ with a trivialization 
$E_{U}\approx U\times E_{\xi}$ such that
$\Pi_{V} |_{E_{U}}$ is isomorphic to the direct product of the
zero Poisson structure and $\Pi_{V} |_{E_{\xi}}$ by a
diffeomorphism $\mathfrak{End}_{B}^{0}(E_{U})$.
\end{proposition}

Note that the data $(\Gamma,F)$ are not preserved under the action
of $\mathfrak{Diff}_{B}^{0}(E)$ in general (compare with Proposition~2.9).

Resuming, we conclude that, according to Theorem~3.2, the
problem on the semilocal equivalence of Poisson structures in
the class $\mathfrak{Coup}_{B}(E)$ is reduced to the study of
coupling Poisson tensors with fixed vertical part. 

\subsection{Sufficient criteria}
Denote by $\mathfrak{Coup}_{B}(E,\Upsilon)$ the set of germs at $B$
of all coupling Poisson tensors $\Pi$ on $E$ with fixed vertical
part, $\Pi_{V}=\Upsilon$, which satisfy (3.1). Then,
\begin{equation}
\rank\Upsilon=0\qquad\text{on}\quad B.
\tag{3.13}
\end{equation}

Denote by $\mathfrak{F}_{B}$ the space of all germs at $B$ of all
smooth functions on vanishing at $B$. Introduce also the subspace 
$\Casim_{B}(E,\Upsilon)$ of $C_{B}^{\infty}(E)$ consisting of
all Casimir functions of $\Upsilon$ vanishing at $B$. Then, we
have a natural identification 
$\Casim_{B}(E,\Upsilon)/\pi^{\ast}C^{\infty}(B)$. The
corresponding space of germs at $B$ of functions in 
$\Casim_{B}(E,\Upsilon)$ will be denoted by $\mathfrak{Casim}_{B}$.

Let $\Pi\in\mathfrak{Coup}_{B}(E,\Upsilon)$ be a coupling Poisson tensor
associated with geometric data $(\Upsilon,\Gamma,F)$. By (3.3)
it is clear that the restriction of the covariant exterior
derivative $\partial^{\Gamma}$ (2.2) to 
$\Omega^{k}(B)\otimes\mathfrak{Casim}_{B}$ is well defined. 
Thus, we have the coboundary operator
\begin{equation}
\partial_{0}^{\Gamma}:\Omega^{k}(B)\otimes\mathfrak{Casim}_{B}
\rightarrow
\Omega^{k+1}(B)\otimes\mathfrak{Casim}_{B}.
\tag{3.14}
\end{equation}

Let $\tilde{\Pi}\in\mathfrak{Coup}_{B}(E,\Upsilon)$ be another
coupling Poisson tensor associated with geometric data
$(\Upsilon,\tilde{\Gamma},\tilde{F})$. 
Recall that we write $\tilde{\Pi}\sim_{Q}\Pi$ if there exists
$Q\in\Omega^{1}(B)\otimes\mathfrak{F}_{B}$ satisfying (2.60). 
This implies $\partial_{0}^{\Gamma}=\partial_{0}^{\tilde{\Gamma}}$. 
By (3.4) we have
\begin{equation}
\pi^{\ast}(\tilde{F}-F) |_{B}=0.
\tag{3.15}
\end{equation}
If $\tilde{\Pi}\sim_{Q}\Pi$, then formula (2.87) defines the
Casimir $2$-cocycle associated to the pair $(\tilde{\Pi},\Pi)$:
\begin{equation}
C=C_{\tilde{\Pi},\Pi}\in\Omega^{2}(B)\otimes\mathfrak{Casim}_{B},
\qquad \partial_{0}^{\Gamma}C=0.
\tag{3.16}
\end{equation}
Recall that $Q$ is uniquely determined by (2.84) up to the
transformation $Q\mapsto Q+K$, for arbitrary
$K\in\Omega^{1}(B)\otimes\mathfrak{Casim}_{B}$. 
Then $C$ is transformed by the rule 
$C\mapsto C+\partial_{0}^{\Gamma}K$
preserving the $\partial_{0}^{\Gamma}$-cohomology class~$[C]$.

We start with the following result which was stated in~\cite{Vo1}.

\begin{theorem}
Let $\Pi,\tilde{\Pi}\in\mathfrak{Coup}_{B}(E,\Upsilon)$ be two
coupling Poisson tensors such that
\begin{equation}
\tilde{\Pi}\sim_{Q}\Pi
\tag{3.17}
\end{equation}
and
\begin{equation}
[ C_{\tilde{\Pi},\Pi}]=0.
\tag{3.18}
\end{equation}
Then, $\Pi$ and $\tilde{\Pi}$ are isomorphic by a diffeomorphism
$\phi \in\mathfrak{Diff}_{B}^{0}(E)$,
\begin{equation}
\phi^{\ast}\tilde{\Pi}=\Pi.
\tag{3.19}
\end{equation}
\end{theorem}

The proof is divided in few steps. According to the results of
Section~2, it suffices to find a Poisson homotopy.

\textit{Step} I. Let $(\Upsilon,\Gamma,F)$ and 
$(\tilde{\Upsilon},\tilde{\Gamma},\tilde{F})$ be the geometric
data of $\Pi$ and $\tilde{\Pi}$, respectively. 
By conditions (3.17), (3.18), one can choose 
$Q\in\Omega^{1}(B)\otimes\mathfrak{F}_{B}$ such that
\begin{align}
\hor^{\tilde{\Gamma}}(u)
&=\hor^{\Gamma}(u)+\Upsilon^{\natural}dQ(u),
\tag{3.20}
\\
\tilde{F}
&=F-(\partial^{\Gamma}Q+\frac{1}{2}\{Q\wedge Q\}_{\Upsilon}).
\tag{3.21}
\end{align}
Define
\begin{equation}
\Gamma^{t}\overset{\text{def}}{=}\Gamma-t\Upsilon^{\sharp}dQ
\tag{3.22}
\end{equation}
and
\begin{equation}
F^{t}\overset{\text{def}}{=}
F-(t\partial^{\Gamma}Q+\frac{t^{2}}{2}\{Q\wedge Q\}_{\Upsilon}).
\tag{3.23}
\end{equation}

\begin{lemma}
For any $t\in[0,1]$, the data $(\Upsilon,\Gamma^{t},F^{t})$ satisfy
relations {\rm(2.23)--(2.26)}.
\end{lemma}

\begin{proof}
It follows from (2.63) that
$$
\partial^{\Gamma^{t}}\Theta
=\partial^{\Gamma}\Theta+t\{Q\wedge\Theta\}_{\Upsilon}.
$$
From here and by using (2.80)--(2.82), we get
\begin{align*}
\partial^{\Gamma^{t}}F^{t}  
& =\partial^{\Gamma}F_{t}+t\{Q\wedge F_{t}\}_{\Upsilon}\\
& =t[-(\partial^{\Gamma})^{2}Q+\{Q\wedge F\}_{\Upsilon}]
-t^{2}\bigg[\frac{1}{2}\partial^{\Gamma}\{Q\wedge Q\}_{\Upsilon}
+\{Q\wedge \partial^{\Gamma}Q\}_{\Upsilon}\bigg]\\
&\qquad 
-\frac{t^{3}}{2}\{Q\wedge\{Q\wedge Q\}_{\Pi_{V}}\}_{\Pi_{V}}=0.
\end{align*}
Next,
$$
L_{^{\hor^{\Gamma^{t}}(u)}}\Upsilon
=L_{\hor^{\Gamma}(u)}\Upsilon+tL_{\Upsilon^{\natural}dQ(u)}\Upsilon=0.
$$
Finally, using (2.71), it is easy to verify that the curvature
$2$-form of $\Gamma_{t}$ satisfies (2.26) precisely for $F^{t}$
in~(3.23). 
\end{proof}

Now, we remark that $F^{t}$ satisfies the nondegeneracy
condition in a neighborhood of $B$ in $E$ for all $t\in[0,1]$. 
Thus, for every $t$, the triple
$(\Upsilon,\Gamma^{t},F^{t})$ defines the coupling Poisson
tensor $\Pi_{t} \in\mathfrak{Coup}_{B}(E,\Upsilon)$.

\textit{Step} II. The family $\{\Pi_{t}\}$ satisfies (2.60) and
(2.61) for $Q^{t}=tQ$. Moreover, $\Pi_{1}=\tilde{\Pi}$ and
$\Pi_{0}=\Pi$. In other words, 
$\Pi_{t}\approx_{Q^{t}}\Pi_{0}$ and $[C_{\Pi_{t}\Pi_{0}}]=0 $
for every~$t$. 
Define the time-dependent vector field 
$X_{t}=X_{t}^{i}(\xi,x)\hor_{i}^{t}$ by
\begin{equation}
\mathbf{i}_{Xt}F^{t}=Q.
\tag{3.24}
\end{equation}
Then,
\begin{equation}
X_{t} |_{B}=0
\tag{3.25}
\end{equation}
and hence for every $t\in[0,1]$, the flow $\Phi_{t}$ of $X_{t}$
is well defined in a neighborhood $N_{\flow}$ of $B$ in~$E$.
According to Theorem~2.18, the time-$1$ flow $\phi=\Phi_{1}$
gives an isomorphism in (3.19). This completes the proof
of the theorem.
\qed

\begin{remark}
It follows from (2.41), that the flow $\Phi_{t}$ does not
preserve the fiber-tangent Poisson structure $\Upsilon$, in general.
\end{remark}

\textbf{Deformations of the coupling form}. 
Let $\Pi,\tilde{\Pi}$ be two
coupling Poisson structures satisfying (3.17). Then the Casimir
$2$-cocycle $C=C_{\tilde{\Pi},\Pi}$ vanishes on~$B$,
$$
\pi^{\ast}C |_{B}=0.
$$
Thus, the $2$-form $F+C\in\Omega^{2}(B)\otimes\mathfrak{F}_{B}$
satisfies (3.4) and the nondegeneracy condition in a
neighborhood of $B$ in $E$. Moreover, it is easy to see that the
triple $(\Upsilon,\Gamma,F+C)$ satisfies the relations
(2.23)--(2.26). 

\begin{proposition}
Every pair $\Pi,\tilde{\Pi}\in\mathfrak{Coup}_{B}(E,\Upsilon)$
such that $\tilde{\Pi}\sim_{Q}\Pi$ induces the coupling Poisson
tensor $\Pi_{C} \in\mathfrak{Coup}_{B}(E,\Upsilon)$ associated with
the geometric data $(\Upsilon,\Gamma,F+C)$.
\end{proposition}

We will call $F+C$ the \textit{deformed coupling} form and
$\Pi_{C}$ the \textit{deformed coupling Poisson tensor}. 
It is clear that the relative Casimir $2$-cocycle corresponding
to $\tilde{\Pi}$ and $\Pi_{C}$ is zero.
Applying Theorem~3.4 to $\tilde{\Pi}$ and $\Pi_{C}$, we get the
following result. 

\begin{theorem}
Let $\Pi,\tilde{\Pi}\in\mathfrak{Coup}_{B}(E,\Upsilon)$ be two
coupling Poisson tensors such that $\tilde{\Pi}\sim_{Q}\Pi$. 
Let $C=C_{\tilde{\Pi}\Pi}$ be the Casimir $2$-cocycle and
$\Pi_{C}$ the deformed coupling Poisson tensor 
associated with geometric data $(\Upsilon,\Gamma,F+C)$. 
Then, $\tilde{\Pi}$ and $\Pi_{C}$ are isomorphic under the
time-$1$ flow of the time-dependent vector field $X_{t}$ defined
by 
\begin{equation}
\mathbf{i}_{Xt}(F^{t}+C)=Q.
\tag{3.26}
\end{equation}
\end{theorem}

\textbf{Equivalent vertical parts}. 
Now suppose we are given a second fiber-tangent Poisson
structure $\tilde{\Upsilon}$ vanishing at $B$. Consider 
the corresponding class of coupling Poisson tensors 
$\mathfrak{Coup}_{B} (E,\tilde{\Upsilon})$. 
Let $\Pi\in\mathfrak{Coup}_{B}(E,\Upsilon)$ and 
$\tilde {\Pi}\in\mathfrak{Coup}_{B}(E,\tilde{\Upsilon})$ be two
coupling Poisson tensors. 
Then, $(B,\omega)$ is a common symplectic leaf of $\Pi$ and
$\tilde{\Pi}$. 
Assume that $\Upsilon$ and $\tilde{\Upsilon}$ are \textit{equivalent},
\begin{equation}
g^{\ast}\tilde{\Upsilon}=\Upsilon
\tag{3.27}
\end{equation}
for some $g\in\mathfrak{End}_{B}^{0}(E)$. 
Then, by Proposition~2.9,
$g^{\ast} \tilde{\Pi}$ is a coupling Poisson tensor associated
with geometric data 
$(\Upsilon,g^{\ast}\tilde{\Gamma},g^{\ast}\tilde{F})$. It is
clear that $g^{\ast}\tilde{\Pi}\in\mathfrak{Coup}_{B}(E,\Upsilon)$.
Assume that 
\begin{equation}
g^{\ast}\tilde{\Pi}\sim_{Q}\Pi
\tag{3.28}
\end{equation}
for a certain $Q\in\Omega^{1}(B)\otimes\mathfrak{F}_{B}$. Then, one
can define the Casimir $2$-cocycle
$C_{g}=C_{g^{\ast}\tilde{\Pi},\Pi}$ and the deformed 
coupling tensor $\Pi_{C_{g}}$ associated to
$g^{\ast}\tilde{\Pi}$ and $\Pi$. 

As a consequence of Theorem 3.8 we get the following result.

\begin{proposition}
Under assumptions {\rm(3.27), (3.28)}, the coupling Poisson
tensors $\tilde{\Pi}$ and $\Pi_{C_{g}}$ are equivalent.
\end{proposition}

Let us look at the dependence of the Casimir $2$-cocycle $C_{g}$
on the choice of $g$ in (3.27).

Let
$\mathfrak{Aut}_{B}^{0}(E,\Upsilon)\subset\mathfrak{Diff}_{B}^{0}(E)$ be
the connected component of the group of germs at $B$ of Poisson
automorphisms of $\Upsilon$. Denote by
$\mathfrak{Inn}_{B}(E,\Upsilon)$ the group of germs at $B$ 
of inner automorphisms of $\Upsilon$ which are identical on~$B$.

The diffeomorphism $g$ is uniquely determined by (3.27) up to
the multiplication by an element
$\varphi\in\mathfrak{Aut}_{B}^{0}(E,\Upsilon)$. 
Assume that the fiber-tangent Poisson structure $\Upsilon$
possesses the property
\begin{equation}
\mathfrak{Aut}_{B}^{0}(E,\Upsilon)=\mathfrak{Inn}_{B}((E,\Upsilon).
\tag{3.29}
\end{equation}

\begin{proposition}
The $\partial_{0}^{\Gamma}$-cohomology class of the Casimir
$2$-cocycle $C_{g}=C_{g^{\ast}\tilde{\Pi},\Pi}$ is independent
of the choice of $g$ in {\rm(3.27)}. 
\end{proposition}

This assertion follows from Proposition~2.24.

Thus, under assumption (3.29) and the equivalence of $\Upsilon$
and $\tilde{\Upsilon}$, the cohomology class
$[C_{g^{\ast}\tilde{\Pi},\Pi}]$ is an invariant of two coupling
Poisson tensors $\Pi\in\mathfrak{Coup}_{B} (E,\Upsilon)$ and 
$\tilde{\Pi}\in\mathfrak{Coup}_{B}(E,\tilde{\Upsilon})$.

\subsection{Cohomological obstructions}
Here we show that the cohomology class of the relative Casimir
$2$-cocycle can be viewed as an obstruction to the Poisson
equivalence in the following special case.

Consider again the set $\mathfrak{Coup}_{B}(E,\Upsilon)$ of germs at
$B$ of all coupling Poisson tensors $\Pi$ on $E$ with fixed
vertical part $\Upsilon$ and satisfying~(3.1).

For every $\xi\in B$, the fiber $E_{\xi}$ carries the Poisson
structure  $\Upsilon_{\xi}=\Upsilon |_{\xi}$ vanishing at the
origin $0$, $\Upsilon_{\xi}=0$ at $0$. 
Since $\Upsilon$ is locally trivial, the germ of
$\Upsilon_{\xi}$ at $0$ is independent of $\xi$ up to an
isomorphism. Fix $\xi^{0}\in B$ and make the following
assumption on $(E_{\xi^{0}},\Upsilon_{\xi^{0}})$. 
Consider the Lichnerowicz complex 
$\mathcal{D}_{k}=[\Upsilon_{\xi^{0}},\cdot]:\,
\chi^{k}(E_{\xi^{0}})\rightarrow\chi^{k+1}(E_{\xi^{0}})$ 
induced by the Schouten bracket on $E_{\xi^{0}}$
\cite{KM,Va1}. In particular,
\begin{align*}
\mathcal{D}_{0}f&=[\Upsilon_{\xi^{0}},f]\equiv\Upsilon_{\xi^{0}}^{\natural}df,
\\
\mathcal{D}_{1}w&=[\Upsilon_{\xi^{0}},w]=L_{w}\Upsilon_{\xi^{0}}
\end{align*}
for $f\in\chi^{0}(E_{\xi^{0}})=C^{\infty}(E_{\xi^{0}})$ and for
$w\in\chi^{1}(E_{\xi^{0}})=\mathcal{X}^{1}(E_{\xi^{0}})$. 
We assume that there exists a neighborhood $O$ of $0$ in
$E_{\xi^{0}}$ and linear operators  
$$
\mathcal{H}_{0}:\,\chi^{1}(\bar{O})\rightarrow\chi^{0}(\bar{O})
\qquad\text{and}\qquad
\mathcal{H}_{1}:\,\chi^{2}(\bar{O})\rightarrow\chi^{1}(\bar{O})
$$
satisfying the homotopy condition
\begin{equation}
\mathcal{D}_{0}\circ\mathcal{H}_{0}
+\mathcal{H}_{1}\circ\mathcal{D}_{1}=\id.
\tag{3.30}
\end{equation}
This assumption guarantees that every smooth family of Poisson
vector fields on $\bar{O}$ is a family of Hamiltonian vector
fields, where the parameter-dependent Hamiltonian smoothly
varies with a parameter. Then, by the standard partition of
unity argument and the local triviality of $\Upsilon$ we 
deduce that property (2.92) holds in a neighborhood of $B$ in
$E$. In particular, for any 
$\Pi,\tilde{\Pi}\in\mathfrak{Coup}_{B}(E,\Upsilon)$ there
exists $Q\in\Omega^{1}(B)\otimes\mathfrak{F}_{B}$ such that 
$\tilde{\Pi}\sim_{Q}\Pi$ and the Casimir $2$-cocycle 
$C_{\tilde{\Pi},\Pi}$ is well defined with respect to the
intrinsic coboundary operator $\partial_{0} $ be associated 
with $\Upsilon$. Remark also that the above assumption implies
(3.29). 

Now we assume that we have the set
$\mathfrak{Coup}_{B}(E,\tilde{\Upsilon})$ corresponding to another
fiber-tangent Poisson structure $\tilde{\Upsilon}$. 
If $\Upsilon$ and $\tilde{\Upsilon}$ are equivalent, then by
Proposition~2.24 for any $\Pi\in\mathfrak{Coup}_{B}(E,\Upsilon)$ 
and $\tilde{\Pi}\in\mathfrak{Coup}_{B}(E,\tilde{\Upsilon})$ 
the $\partial_{0}$-cohomology class
$[C_{g^{\ast}\tilde{\Pi},\Pi}]$ is independent of the choice of
$g$ in (3.27). 
Then we can put
\begin{equation}
[ C_{\tilde{\Pi},\Pi}]\overset{\text{def}}{=}
[C_{g^{\ast}\tilde{\Pi},\Pi}].
\tag{3.31}
\end{equation}
We arrive at the main observation.

\begin{theorem}
The coupling Poisson tensors 
$$
\Pi\in\mathfrak{Coup}_{B}(E,\Upsilon)\qquad\text{and}\qquad
\tilde {\Pi}\in\mathfrak{Coup}_{B}(E,\tilde{\Upsilon})
$$ 
are isomorphic by a $\phi \in\mathfrak{Diff}_{B}^{0}(E)$ if and only if

{\rm(a)} the germs at $B$ of $\Upsilon$ and $\tilde{\Upsilon}$
are equivalent, 

\noindent
and

{\rm(b)} the $\partial_{0}$-cohomology class of $(\tilde{\Pi},\Pi)$ 
is trivial,
$$
[ C_{\tilde{\Pi},\Pi}]=0.
$$
\end{theorem}

The sufficiency follows from Proposition 3.9. The necessity part is a
consequence of Theorem~2.19.

\begin{example}
If $E_{\xi^{0}}=\mathfrak{g}^{\ast}$ is the dual of a semisimple Lie
algebra $\mathfrak{g}$ of compact type and $\Upsilon_{\xi^{0}}$ is
the Lie--Poisson bracket on $\mathfrak{g}^{\ast}$, then, as shown 
in~\cite{Co2}, for a closed ball centered at $0$ there exist 
homotopy operators in (3.30). 
\end{example}

\section{Linearizability and normal forms}

Here we give some results on the linearization of a Poisson
structure at a given (singular) symplectic leaf. 
First, we recall a coordinate derivation of
the linearized Poisson structure of the leaf~\cite{Vo2}. An
invariant definition related to the notion of the Lie algebroid
of a symplectic leaf~\cite{We2} can be found in \cite{Vo1,Va2}.

\subsection{Linearized Poisson structures}
Let $(M,\Psi)$ be a Poisson manifold and $(B,\omega)$ be a
symplectic leaf. We will assume that $B$ is an embedded
submanifold of $M$. Consider the normal bundle 
$\pi:\,E\rightarrow B$ to $B$,
$$
E=T_{B}M/TB.
$$
As usual, to study Poisson geometry around the leaf $B$, we can
move from $M$ to the fibered space (the total space $E$ of
$\pi$) by means of an exponential map. Consider the decomposition
\begin{equation}
T_{B}E=TB\oplus E
\tag{4.1}
\end{equation}
and denote by $\tau_{b}:\,T_{b}E\rightarrow E$ the projection
along $TB$. We say that a diffeomorphism 
$\mathbf{f}:\,E\rightarrow M$ from the total space onto a
neighborhood of the leaf $B$ in $M$ is an 
\textit{exponential map} if $\mathbf{f} |_{B}=\id_{B}$ and
$$
\boldsymbol{\nu}_{b}\circ d_{b}\mathbf{f}=\tau_{b}
$$
for every $b\in B$. 
Here $\boldsymbol{\nu}_{b}:\,T_{b}M\rightarrow E_{b}$ is the
natural projection. An exponential map exists always because of
the tubular neighborhood theorem~\cite{LMr}.

Pick an exponential map $\mathbf{f}$ and consider the pull-back
$\Pi=\mathbf{f}^{\ast}\Psi$ via $\mathbf{f}$. 
It is clear that $(B,\omega)$ (as the zero section) is a
symplectic leaf of the Poisson tensor $\Pi$. Then by 
Proposition~3.1, $\Pi$ is coupling Poisson tensor on $E$
associated with geometric data $(\Pi_{V},\Gamma,F)$ satisfying
conditions (3.2)--(3.4). 

Let $(\xi,x)=(\xi^{i},x^{\sigma})$ be a (local) coordinate
system on $E$, where $(\xi^{i})$ are coordinates on $B$ and 
$(x^{\sigma})$ are affine coordinates along the fibers, so that
locally, $B=\{x=0\}$. Then, 
\begin{equation}
\Pi 
=\Pi_{H}+\Pi_{V}
=-\frac{1}{2}F^{ij}(\xi,x)\hor_{i}\wedge\hor_{j}
+\frac{1}{2}\Pi_{V}^{\alpha\beta}(\xi,x)
\frac{\partial}{\partial x^{\alpha}}\wedge
\frac{\partial}{\partial x^{\beta}},
\tag{4.2}
\end{equation}
where $\hor_{i}=\frac{\partial}{\partial\xi^{i}}
-\Gamma_{i}^{\nu}(\xi,x)\frac{\partial}{\partial x^{\nu}}$. 
In coordinates, conditions (3.2)--(3.4) read
\begin{align}
\Pi_{V}^{\alpha\beta}(\xi,0)&=0,
\tag{4.3}
\\
\Gamma_{i}^{\nu}(\xi,0)&=0
\tag{4.4}
\\
F^{ij}(\xi,0)&=\omega^{ij}(\xi),
\tag{4.5}
\end{align}

Here
$$
\omega=\frac{1}{2}\omega_{ij}(\xi)d\xi^{i}\wedge d\xi^{j}.
$$

\textbf{Linearization of $\Pi_{V}$.} 
The Taylor expansion at $x=0$ for $\Pi_{V}$ gives
$$
\Pi_{V}=\Lambda+O_{2},
$$
where
\begin{equation}
\Lambda=\frac{1}{2}\lambda_{\gamma}^{\alpha\beta}(\xi)x^{\gamma}
\frac{\partial}{\partial x^{\alpha}}\wedge
\frac{\partial}{\partial x^{\beta}}
\tag{4.6}
\end{equation}
and
$$
\lambda_{\gamma}^{\alpha\beta}(\xi)\overset{\text{def}}{=}
\frac{\partial}{\partial x^{\sigma}}\Pi_{\gamma}^{\alpha\beta}(\xi,0).
$$
It is clear that $\Lambda$ is a global bivector field on $E$.
Linearizing the Jacobi identity for $\Pi_{V}$ leads to
$[\Lambda,\Lambda]=0$. Moreover, $\Lambda$ is independent of the
choice of an exponential map. Thus, $\Lambda$ is an intrinsic
fiber-tangent Poisson structure on $E$ of the leaf $B$. The 
restriction of $\Lambda$ to each fiber $E_{\xi}$ gives the
Lie--Poisson structure called the 
\textit{linearized transverse Poisson structure} of
$B$~\cite{We1}. This bundle of Lie--Poisson manifolds is locally
trivial with typical fiber $\mathfrak{g}^{\ast}$, where $\mathfrak{g}$
is the \textit{transverse Lie algebra} of the leaf~$B$.

\textbf{Linearization of $\Gamma$.} 
Consider the connection $\Gamma$,
$$
\Gamma=\big(dx^{\sigma}+\Gamma_{i}^{\sigma}(\xi,x)d\xi\big)\otimes
\frac{\partial}{\partial x^{\sigma}}.
$$
By (4.4) we have
$$
\Gamma_{i}^{\sigma}(\xi,x)=\theta_{i\nu}^{\sigma}(\xi)x^{\nu}+O(x^{2}).
$$
Here the coefficients $\theta_{i\nu}^{\sigma}(\xi)$ are defined in terms of
the Poisson bracket of $\Pi$ as follows
\begin{align*}
\theta_{i\nu}^{\sigma}(\xi)  
=\frac{\partial\Gamma_{i}^{\sigma}(\xi,0)}{\partial x^{\nu}}
=\omega_{ij}(\xi)
\bigg[\frac{\partial}{\partial x^{\nu}}\{\xi^{j},x^{\sigma}\}_{\Pi}\bigg]  
\bigg |_{x=0}.
\end{align*}
Then,
\begin{equation}
\Gamma^{(1)}=(dx^{\sigma}+\theta_{i\nu}^{\sigma}
(\xi)x^{\nu}d\xi^{i})\otimes\frac{\partial}{\partial x^{\sigma}}
\tag{4.7}
\end{equation}
is an \textit{homogeneous} Ehresmann connection on $E$, in the
sense that the Lie derivative along the horizontal lift
$\hor^{\Gamma^{(1)}}(u)$ preserves the space of fiberwise linear
functions on $E$. The linearization of (2.24) leads to
\begin{equation}
L_{\hor^{\Gamma^{(1)}}(u)}\Lambda=0.
\tag{4.8}
\end{equation}
One can associate to $\Gamma^{(1)}$ the \textit{linear
connection} $\nabla$ on $E$ determined by the matrix-valued
connection form $\theta=(\theta_{i\nu}^{\sigma}(\xi)d\xi^{i})$. 
In terms of $\nabla$, condition (4.9) means that
the bivector field $\Lambda$ is 
$\nabla$-\textit{covariantly constant},
$\nabla\Lambda=0$. A linear connection satisfying this property
is called a \textit{Lie--Poisson connection} on~$(E,\Lambda)$.

\textbf{Linearization of $F$.} 
In coordinates
$$
\pi^{\ast}F=\frac{1}{2}F_{ij}(\xi,x)d\xi^{i}\wedge d\xi^{j}.
$$
Taking into account (4.5), we get
$$
F_{ij}(\xi,x)
=\omega_{ij}(\xi)-\mathcal{R}_{ij\sigma}(\xi)x^{\sigma}+O(x^{2}),
$$
where
$$
\mathcal{R}_{ij\sigma}(\xi)=-\omega_{ii^{\prime}}(\xi)
\bigg[\frac{\partial}{\partial x^{\sigma}}
\{\xi^{i^{\prime}},\xi^{j^{\prime}}\}_{\Pi}\bigg]
\bigg |_{x=0}\omega_{j^{\prime}j}(\xi).
$$
We define $F^{(1)}\in\Omega^{1}(B)\otimes C^{\infty}(E)$ by
\begin{equation}
\pi^{\ast}F^{(1)}=\frac{1}{2}F_{ij}^{(1)}(\xi,x)d\xi^{i}\wedge 
d\xi^{j}
\tag{4.9}
\end{equation}
with
\begin{equation}
F_{ij}^{(1)}(\xi,x)\overset{\text{def}}{=}
\omega_{ij}(\xi)-\mathcal{R}_{ij\sigma}(\xi)x^{\sigma}.
\tag{4.10}
\end{equation}
It follows from (2.25) and (2.26) that
\begin{equation}
\partial^{\Gamma^{(1)}}F^{(1)}=0
\tag{4.11}
\end{equation}
and
\begin{equation}
\Curv^{\Gamma^{(1)}}(u_{1},u_{2})
=\Lambda^{\sharp}dF^{(1)}(u_{1},u_{2}).
\tag{4.12}
\end{equation}
Now we can define the following bivector field
\begin{align} 
& \Pi^{(1)}\overset{\text{def}}{=}\Pi_{H}^{(1)}+\Lambda
\tag{4.13}\\
&\qquad 
=-\frac{1}{2}(F^{(1)})^{ij}(\xi,x)\hor_{i}\wedge
\hor_{j}+\frac{1}{2}\lambda_{\gamma}^{\alpha\beta}(\xi)x^{\gamma}
\frac{\partial}{\partial x^{\alpha}}\wedge
\frac{\partial}{\partial x^{\beta}}
\nonumber
\end{align}
where $(F^{(1)})^{is}F_{sj}^{(1)}=\delta_{j}^{i}$ and
$$
\hor_{i}=\hor_{i}^{\Gamma^{(1)}}=\frac{\partial}{\partial\xi^{i}}
-\theta_{i\nu}^{\sigma}(\xi)x^{\nu}
\frac{\partial}{\partial x^{\sigma}}.
$$
The bivector $\Pi^{(1)}$ is well defined on the following
neighborhood of~$B$: 
$$
N^{(1)}=\{(\xi,x)\in E\mid\det[F_{ij}^{(1)}(\xi,x)]\neq0\}.
$$

Taking into account (4.8), (4.11), and (4.12), by
Proposition~2.4 we get the following fact.

\begin{proposition}
Bivector field $\Pi^{(1)}$ in {\rm(4.13)} is the coupling
Poisson tensor associated to the geometric data 
$(\Lambda,\Gamma^{(1)},F^{(1)})$,
$$
\Pi^{(1)}\in\mathfrak{Coup}_{B}(E,\Lambda).
$$
\end{proposition}

In particular, $(B,\omega)$ is a symplectic leaf of $\Pi^{(1)}$ and
$$
\Pi=\Pi^{(1)}+O_{2},
$$
that is, $\Pi^{(1)}$ is a \textit{first approximation} to $\Pi$ at~$B$.

Let $E^{\ast}$ be the dual bundle of $E$ called the co-normal
bundle of the leaf $B$. Then, $E^{\ast}$ is a bundle of Lie
algebras with typical fiber $\mathfrak{g}$. 
Using (4.10), one can introduce the vector-valued $2$-form
$\mathcal{R}\in\Omega^{2}(B,\Sect(E^{\ast}))$, locally given by
\begin{equation}
\mathcal{R}=\frac{1}{2}\mathcal{R}_{ij\sigma}(\xi)d\xi^{i}\wedge 
d\xi^{j}\otimes e^{\sigma}(\xi).
\tag{4.14}
\end{equation}
Here $\{e^{\sigma}\}$ is the basis of local sections of $E^{\ast}$
corresponding to the coordinates $x^{\sigma}$ in (4.10). Let
$\Curv^{\nabla}\in\Omega^{2}(B,\End(E^{\ast}))$
be the curvature form of the linear connection $\nabla$,
$\Curv^{\nabla}(u_{1},u_{2})
=[\nabla_{u_{1}},\nabla_{u_{2}}]-\nabla_{[u_{1},u_{2}]}$. 
Then, we have
$$
\Curv^{\Gamma^{(1)}}(u_{1},u_{2})
=\big\langle\Curv^{\nabla}(u_{1},u_{2})x,\frac{\partial}{\partial x}
\big\rangle.
$$
It follows from here that equality (4.12) can be rewritten in the
form
\begin{equation}
\Curv^{\nabla}=-\ad{\ast}\circ\mathcal{R},
\tag{4.15}
\end{equation}
where $\ad^{\ast}$ is the coadjoint operator on the fibers of
$E$. Thus, the coupling Poisson tensor $\Pi^{(1)}$ is determined
by the data $(\Lambda,\nabla,\mathcal{R})$.

\begin{remark}
One can define \cite{Vo1} the Poisson tensor $\Pi^{(1)}$starting
with the transitive Lie algebroid $T_{B}^{\ast}E$ 
of the symplectic leaf~\cite{We2}. 
Then, $\Pi^{(1)}$ is completely determined by the choice of a
connection of the Lie algebroid \cite{Ku,Mz}. 
The $2$-form $\mathcal{R}$ is just the curvature of this connection.
\end{remark}

\textbf{Varying the exponential map.} 
Suppose we are given two exponential
maps $\mathbf{f}:\,E\rightarrow M$ and
$\tilde{\mathbf{f}}:\,E\rightarrow M$. 
Consider the corresponding coupling Poisson tensors
$\Pi=\mathbf{f}^{\ast}\Psi$ and 
$\tilde{\Pi}=\tilde{\mathbf{f}}^{\ast}\Psi$. 
Then $\Pi=\phi^{\ast}\tilde{\Pi}$, where 
$\phi=\tilde{\mathbf{f}}^{-1}\circ\mathbf{f}
\in\mathfrak{Diff}_{B}^{0}(E)$. 
Indeed, by the definition of the exponential map,
we have $\phi |_{B}=\id_{B}$ and
$$
d_{\xi}\phi |_{E_{\xi}}=\id_{E_{\xi}}
$$
for every $\xi\in B$. Recall that we have the decomposition
(4.1). Using the dilation $m_{t}:\,E_{\xi}\rightarrow E_{\xi}$ 
by the factor $t\in[0,1]$, $m_{t}(x)=t\cdot x$, 
one can define 
$\Phi_{t}=m_{t}^{-1}\circ\phi\circ m_{t}$. It
follows that $\{\Phi_{t}\}_{t\in[0,1]}$ is a smooth family of
diffeomorphisms with properties: $\Phi_{0}=\id$, 
$\Phi_{1}=\phi$ and $\Phi_{t} |_{B}=\id_{B}$. Then, $\Phi_{t}$
is the flow of the time-dependent vector field 
$Z_{t}=\frac{d\Phi_{t}}{dt}\circ\Phi_{t}^{-1}$ vanishing at~$B$.

From here and Theorem 3.2 we get the following fact.

\begin{proposition}
The vertical part $\Pi_{V}$ of the coupling Poisson tensor
$\Pi=\mathbf{f}^{\ast}\Psi$ is independent of the choice of an
exponential map $\mathbf{f}$ up to an isomorphism in 
$\mathfrak{End}_{B}^{0}(E)$.
\end{proposition}

Now consider the linearized Poisson structures 
$\Pi^{(1)}$ and $\tilde{\Pi}^{(1)}$ associated to the geometric
data $(\Lambda,\Gamma^{(1)},F^{(1)})$ and 
$(\Lambda,\tilde{\Gamma}^{(1)},\tilde{F}^{(1)})$, respectively.
Consider the decompositions
$$
T_{B}M=TB\oplus\mathcal{L}=TB\oplus\tilde{\mathcal{L}},
$$
where
$$
\mathcal{L}=d_{B}\mathbf{f}(E)\qquad\text{and}\qquad
\tilde{\mathcal{L}}=d_{B}\tilde{\mathbf{f}}(E).
$$
Let $l:\,T_{B}M\rightarrow\tilde{\mathcal{L}}$ be the projection
along $TB$. Pick a basis of (local) sections $\{n_{\sigma}(\xi)\}$ 
of $\mathcal{L}$.
Then, we have
$$
l_{\xi}(n_{\sigma}(\xi))=n_{\sigma}(\xi)+u_{\sigma}(\xi),
$$
where $u_{\sigma}(\xi)\in T_{\xi}B$. Define 
$Q^{(1)}\in\Omega^{1}(B)\otimes C^{\infty}(E)$ by
\begin{equation}
Q^{(1)}\overset{\text{def}}{=}-(\mathbf{i}_{u_{\nu}}\omega)x^{\nu}.
\tag{4.16}
\end{equation}
Let $\partial^{\Gamma^{(1)}}$ be the exterior covariant
derivative associated to $\Gamma^{(1)}$.

\begin{proposition}
The linear connections $\Gamma^{(1)}$ and $\tilde{\Gamma}^{(1)}$
are related by
\begin{equation}
\hor^{\tilde{\Gamma}^{(1)}}(u)=\hor^{\Gamma^{(1)}}(u)
+\Lambda^{\sharp}dQ^{(1)}(u),
\tag{4.17}
\end{equation}
where $Q^{(1)}$ is given by {\rm(4.16)}. 
The corresponding coupling $2$-forms
$F^{(1)}$ and $\tilde{F}^{(1)}$ satisfy
\begin{equation}
\tilde{F}^{(1)}=F^{(1)}
-\bigg(\partial^{\Gamma^{(1)}}Q^{(1)}
+\frac{1}{2}\{Q^{(1)}\wedge Q^{(1)}\}_{\Lambda}\bigg).
\tag{4.18}
\end{equation}
\end{proposition}

This implies that for 
$\tilde{\Pi}^{(1)},\Pi^{(1)}\in\mathfrak{Coup}_{B}(E,\Lambda)$ 
we have that $\tilde{\Pi}^{(1)}\sim_{Q^{(1)}}\Pi^{(1)}$ and 
the relative Casimir cohomology class is trivial,
$$
[ C_{\tilde{\Pi}^{(1)}\Pi^{(1)}}]=0.
$$
Applying Theorem 3.4 leads to the following consequence.

\begin{corollary}
The germs at $B$ of the coupling Poisson structures $\Pi^{(1)}$
and $\tilde{\Pi}^{(1)}$ are isomorphic by a diffeomorphism in
$\mathfrak{Diff}_{B}^{0}(E)$.
\end{corollary}

Thus, $\Pi^{(1)}$ is independent of the choice of an exponential
map up to isomorphism. We call $\Pi^{(1)}$ the 
\textit{linearized Poisson structure} of
the original Poisson structure $\Psi$ at a given symplectic leaf
$(B,\omega)$. 

Another important consequence of Proposition~4.4 is that the
Poisson structure $\Psi$ with a given symplectic leaf $B$
inherits an \textit{intrinsic coboundary operator},
\begin{equation}
\partial_{0}:\,\Omega^{k}(B)\otimes\mathfrak{Casim}_{B}(E)\rightarrow
\Omega^{k+1}(B)\otimes\mathfrak{Casim}_{B}(E),
\tag{4.19}
\end{equation}
where $\mathfrak{Casim}_{B}(E)$ is the space of germs at $B$ of Casimir
functions of $\Lambda$ vanishing on $B$. We can put 
$\partial_{0}=\partial_{0}^{\Gamma^{(1)}}$, but this definition
is independent of the connection $\Gamma^{(1)}$ because of (4.17).

\subsection{Semilocal linearization problem}
The linearized Poisson structure of the symplectic leaf is a
natural candidate for setting of the semilocal linearization problem.

\begin{definition}
The Poisson structure $\Psi$ is said to be \textit{linearizable} at the
symplectic leaf $(B,\omega)$ if there exists an exponential map
$\mathbf{f}:\, E\rightarrow M$ such that the germs at $B$ of the
pull-back Poisson structure $\Pi=\mathbf{f}^{\ast}\Psi$ and the
linearized Poisson structure $\Pi^{(1)}$ are isomorphic by a
diffeomorphism in $\mathfrak{Diff}_{B}^{0}(E)$. Respectively, 
$\Psi$ is \textit{transversally linearizable} at $B$ if the
germs of $\Pi_{V}$ and $\Lambda$ are isomorphic by a
diffeomorphism in $\mathfrak{End}_{B}^{0}(E)$. 
\end{definition}

It follows from Proposition~4.3 and Corollary~4.5 that this
definition is independent of the choice of the exponential map.
Remark also that this definition agrees with the
zero-dimensional case when $\dim B=0$. 

As a consequence of Theorem 3.2 we get the following result.

\begin{proposition}
The linearizability implies the transversal linearizability.
\end{proposition}

In fact, the linearizability of $\Pi_{V}$ at $B$ is equivalent
to the linearizability of the transverse Poisson structure of a
point in $B$. This statement can be established in various ways 
(see, for example,~\cite{Br}). Our arguments are based on the
local splitting theorem. Let $\mathfrak{g}$ be the transverse Lie
algebra of the symplectic leaf $B$. Recall that the bundle $E$
of Lie--Poisson manifolds associated to the linearized
transverse Poisson structure $\Lambda$ is locally trivial with typical
fiber $\mathfrak{g}^{\ast}$.
In particular, the restriction
$\Lambda_{\xi}=\Lambda |_{E_{\xi}}$ is isomorphic to the
Lie--Poisson structure of $\mathfrak{g}^{\ast}$. Consider 
also the restriction $(\Pi_{V})_{\xi}$ of $\Pi_{V}$ to the fiber
$E_{\xi}$. 

\begin{proposition}
Fix a point $\xi^{0}\in B$. Then the following assertions are
equivalent{\rm:} 

{\rm(a)} $(\Pi_{V})_{\xi^{0}}$ and $\Lambda_{\xi^{0}}$ are
isomorphic by a (local) diffeomorphism $\psi$ such that
$$
\psi(0)=0\qquad\text{and}\qquad d_{0}\psi=\id;
$$

{\rm(b)} the germs at $B$ of $\Pi_{V}$ and $\Lambda$ are
isomorphic by a diffeomorphism in $\mathfrak{End}_{B}^{0}(E)$.
\end{proposition}

\begin{proof}
The first part $(b)\Longrightarrow(a)$ is evident. To prove the
implication $(a)\Longrightarrow(b)$ we consider the following
family of fiber-tangent Poisson structures
$$
\Upsilon_{t}=t\mathfrak{m}_{t}^{\ast}\Pi_{V}
$$
with $\Upsilon_{0}=\Lambda$ and $\Upsilon_{1}=\Pi_{V}$. It
suffices to show that there exists a time-dependent vector field
$Y_{t}$ vanishing at $B$ and satisfying the equation
\begin{equation}
L_{Y_{t}}\Upsilon_{t}+\frac{\partial\Upsilon_{t}}{\partial t}=0.
\tag{4.20}
\end{equation}
First, let us solve Eq.~(4.20) locally. By Proposition~3.3,
in a neighborhood of each point $\xi^{0}\in B$ in $E$
there exists a coordinate system $(\xi^{i},x^{\sigma})$ such that
$$
\Upsilon_{t}=\frac{1}{2t}\Pi_{V}^{\alpha\beta}(tx)
\frac{\partial}{\partial x^{\alpha}}\wedge
\frac{\partial}{\partial x^{\beta}}.
$$
Then, by condition (a) one can choose a local solution of (4.20)
in the form $Y_{t}=Y_{t}(x)\frac{\partial}{\partial x^{\alpha}}$. 
Finally, using a partition of unity on $B$, we can splice
together local solutions to a global one. 
\end{proof}

\textbf{Flat linearized transverse Poisson structures.} 
We say that $\Lambda$ is \textit{flat} if there exists a flat
Lie--Poisson connection $\nabla^{\flat}$ on
$(E,\Lambda)$,  $\Curv^{\nabla^{\flat}}=0$. 
Locally, $\Lambda$ is flat because of the local triviality
property.  Given a flat connection $\nabla^{\flat}$ one can
define the following coupling Poisson tensor
\begin{equation}
\Pi^{\flat} =\Pi_{H}^{\flat}+\Lambda
 =-\frac{1}{2}\omega^{ij}(\xi)\frac{\partial}{\partial\xi^{i}}
\wedge
\frac{\partial}{\partial\xi^{j}}
+\frac{1}{2}\lambda_{\gamma}^{\alpha\beta}x^{\gamma}
\frac{\partial}{\partial x^{\alpha}}\wedge
\frac{\partial}{\partial x^{\beta}}.
\tag{4.21}
\end{equation}
Here $(\xi^{i},x^{\sigma})$ are a coordinate system on $E$ such
that the coordinates $x^{\sigma}$ along the fibers correspond to
a $\nabla^{\flat}$-parallel basis $\{e_{\sigma}\}$ of sections
of $E$, $\nabla^{\flat}e_{\sigma}=0$. Moreover, 
$\lambda_{\gamma}^{\alpha\beta}=\const$ are the structure
constants relative to the dual basis $\{e^{\sigma}\}$. 
The coupling form $F^{\flat}$ of $\Pi^{\flat}$ is given by
$F^{\flat}=\omega\otimes1$. The horizontal part
$\Pi_{H}^{\flat}$ is the lifting of the nondegenerate Poisson 
structure on the leaf $B$ via the connection $\nabla^{\flat}$.

\begin{proposition}
Let $\Pi^{(1)}$ be a linearized Poisson tensor associated with
data $(\Lambda,\nabla,\mathcal{R})$. Assume that the following
conditions hold{\rm:} 

{\rm(a)} $\Lambda$ is flat{\rm;}

{\rm(b)} the center and the first cohomology group of $\mathfrak{g}$
are trivial, 
\begin{align}
Z(\mathfrak{g})&=\{0\},
\tag{4.22}
\\
H^{1}(\mathfrak{g},\mathfrak{g)}&=\{0\}.
\tag{4.23}
\end{align}
Then, $\Pi^{(1)}$ and $\Pi^{\flat}$ are isomorphic by a
diffeomorphism in $\mathfrak{Diff}_{B}^{0}(E)$. If, in addition
to the conditions~{\rm(a)}, and~{\rm(b)}, the holonomy group of
the connection $\nabla^{\flat}$ is trivial, then the normal
bundle $E$ is trivial and $\Pi^{(1)}$ is equivalent to the
direct product Poisson structure on
$E=B\times\mathfrak{g}^{\ast}$. 
\end{proposition}

The proof of this statement is based on Theorem~3.4 and the
observation that conditions (4.23), (4.23) imply properties
(4.17), (4.18) for the connections $\nabla$ and $\nabla^{\flat}$
(for more details, see~\cite{Vo4}). 

Remark that (4.22), (4.23) hold automatically in the case when
$\mathfrak{g}$ is semisimple. 
Some examples, where $\Lambda$ is not flat, 
can be found in~\cite{IKV}.

\subsection{Main results}
Suppose we start with a Poisson manifold $(M,\Psi)$ and an
embedded symplectic leaf $(B,\omega)$. 
Let $(\mathfrak{g},\Lambda,\partial_{0})$ be the intrinsic data
of the leaf consisting of the transverse Lie algebra $\mathfrak{g}$,
linearized transverse Poisson structure $\Lambda$, and the
coboundary operator $\partial_{0}$ in (4.19). 

Assume that $\mathfrak{g}$ is \textit{semisimple of compact type}.
Fix an exponential map $\mathbf{f}$. Consider the pull-back
Poisson structure 
$\Pi=\mathbf{f}^{\ast}\Psi\in\mathfrak{Coup}_{B}(E,\Pi_{V})$ and the
corresponding linearized Poisson structure 
$\Pi^{(1)}\in\mathfrak{Coup}_{B}(E,\Lambda)$. 
Let $(\Pi_{V},\Gamma,F)$ and $(\Lambda,\Gamma^{(1)},F^{(1)})$ 
be the corresponding geometric data.

By the local linearization theorem due to \cite{Co2} and
Proposition~4.8 we deduce the following fact.

\begin{proposition}
$\Psi$ is transversally linearizable at $B$, that is, there exits a
diffeomorphism such that $g\in\mathfrak{End}_{B}^{0}(E)$
\begin{equation}
g^{\ast}\Pi_{V}=\Lambda,
\tag{4.24}
\end{equation}
\end{proposition}

Pick $g$ in (4.24) and consider 
$g^{\ast}\Pi\in\mathfrak{Coup}_{B}(E,\Lambda) $
associated with\break $(\Lambda,g^{\ast}\Gamma,g^{\ast}F)$. 
By assumption, the typical fiber of $(E,\Lambda)$ is the
Lie--Poisson structure of the dual $\mathfrak{g}^{\ast}$ of the
semisimple Lie algebra of compact type. Recall that, according
to~\cite{Co2}, in each closed ball of $0$ in
$\mathfrak{g}^{\ast}$ there exist a homotopy operators in (3.30)
and hence conditions (2.92), (3.29) hold. Thus, there exists
$Q\in\Omega^{1}(B)\otimes\mathfrak{F}_{B}$ such that 
\begin{equation}
\hor^{g^{\ast}\Gamma}(u)=\hor^{\Gamma^{(1)}}(u)+\Lambda^{\sharp}dQ(u).
\tag{4.25}
\end{equation}

Now, one can define the Casimir $2$-cocycle 
$C\in\Omega^{2}(B)\otimes\mathfrak{Casim}_{B}(E)$ by
\begin{equation}
C\overset{\text{def}}{=}g^{\ast}F-F^{(1)}
+\bigg(\partial_{0}Q+\frac{1}{2}\{Q\wedge Q\}_{\Lambda}\bigg).
\tag{4.26}
\end{equation}

\begin{proposition}
The $\partial_{0}$-cohomology class $[C]$ of the Casimir
$2$-cocycle $C$ {\rm(4.26)} is independent of the choice of $g$
in {\rm(4.24)}, $Q$ in {\rm(4.25)} and an exponential map
$\mathbf{f}$. Moreover, 
$$
\pi^{\ast}C=C_{ij}(\xi,x)d\xi^{i}\wedge d\xi^{j},
$$
where the coefficients $C_{ij}$ have zero of the second order at $x=0$,
$$
C_{ij}=O(x^{2}).
$$
\end{proposition}

\begin{proof} 
The independence of $[C]$ of the choice of $g$ follows
straightforwardly from Proposition~2.24. 
Now, let $\Pi=\mathbf{f}^{\ast}\Psi$
and $\tilde{\Pi}=\mathbf{\tilde{f}}^{\ast}\Psi$ be coupling
Poisson tensors corresponding to exponential maps 
$\mathbf{f}$ and $\mathbf{\tilde{f}}$, respectively. 
Without of loss of generality, by Propositions~4.3 and~4.10, 
one can assume that $\Pi_{V}=\tilde{\Pi}_{V}=\Lambda$. 
Since $\phi^{\ast}\tilde{\Pi}=\Pi$, for 
$\phi=\tilde{\mathbf{f}}^{-1}\circ\mathbf{f}
\in\mathfrak{Diff}_{B}^{0}(E)$, 
Theorem~3.11 implies that $[C_{\tilde{\Pi}\Pi}]=0$. 
Moreover, $[C_{\tilde{\Pi}^{(1)}\Pi^{(1)}}]=0$. 
Using property (2.91), we get 
$[C_{\tilde{\Pi}\tilde{\Pi}^{(1)}}]=[C_{\Pi\tilde{\Pi}^{(1)}}]
=[C_{\Pi\Pi^{(1)}}]$.
\end{proof}

We shall call $[C]$ the $\partial_{0}$-\textit{cohomology class}
of the leaf~$B$.

Using the $2$-cocycle $C$, we define the deformed coupling form 
$F^{(1)}+C\in\Omega^{2}(B)\otimes\mathfrak{F}_{B}$. 
Consider the coupling Poisson tensor
$\Pi_{C}^{(1)}\in\mathfrak{Coup}_{B}(E,\Lambda)$ associated to the
geometric data $(\Lambda,\Gamma^{(1)},F^{(1)}+C)$.

By Theorem 3.8, Proposition~3.9, and Proposition~4.10, 
we deduce the normal form theorem.

\begin{theorem}
The Poisson tensors $\Pi$ and $\Pi_{C}^{(1)}$ are isomorphic by a
diffeomorphism $\phi\in\mathfrak{Diff}_{B}^{0}(E)$,
\begin{equation}
\phi^{\ast}\Pi=\Pi_{C}^{(1)}.
\tag{4.27}
\end{equation}
\end{theorem}

It is clear that the germ of the original Poisson structure
$\Psi$ at $B\subset M$ is also isomorphic to $\Pi_{C}^{(1)}$.

\begin{remark}
A similar result to Theorem 4.12 was formulated in~\cite{Br} without
specifying the coupling form.
\end{remark}

Now applying Theorem~3.11, we get the semilocal linearizability
theorem which can be considered as a generalization of the local
linearization theorem~\cite{Co2}.

\begin{theorem}
The Poisson structure $\Psi$ is linearizable at $B$ if and only
if the $\partial_{0}$-cohomology class of the leaf $B$ is zero,
\begin{equation}
[C]=0.
\tag{4.28}
\end{equation}
\end{theorem}

\begin{corollary}
If the second $\partial_{0}$-cohomology space is trivial, 
then every Poisson structure in $\mathfrak{Coup}_{B}(E,\Lambda)$ is
linearizable. 
\end{corollary}

\begin{corollary} 
Let $\Pi^{\prime}\in\mathfrak{Coup}_{B}(E,\Lambda)$ 
be a linearizable coupling Poisson tensor associated with
geometric data 
$(\Lambda,\Gamma^{\prime},F^{\prime})$. 
Let $C^{\prime}\!\in\!\Omega^{2}(B)\otimes\mathfrak{Casim}_{B}(E)$ 
be a $\partial_{0}$-cocycle. 
Consider the deformed coupling Poisson tensor 
$\Pi_{-C^{\prime}}^{\prime}$ associated to the data 
$(\Lambda,\Gamma^{\prime},F^{\prime}-C^{\prime})$. 
If $[C^{\prime}]\neq0$, 
then $\Pi_{-C^{\prime}}^{\prime}$ is nonlinearizable.
\end{corollary}
 
It is of interest to give a geometric interpretation 
of the triviality of the second $\partial_{0}$-cohomology space.

\begin{remark}  
In a linearization conjecture due to~\cite{CrFe}, 
an integrability property of the transitive Lie algebroid 
$T_{B}^{\ast}E$ appears as a sufficient condition 
for the linearizability in the case when~$B$ is compact.
\end{remark}

As a consequence of Proposition 4.9 and Theorem 4.14 we get the
following semilocal analog of the local splitting theorem.

\begin{theorem}
Assume that the linearized transverse Poisson structure
$\Lambda$ is globally trivial, that is, $\Lambda$ admits a flat
Lie--Poisson connection with trivial holonomy. If {\rm(4.28)} 
holds,
then, the germ of $\Psi$ at $B$ is isomorphic to the direct
product Poisson structure on $E=B\times\mathfrak{g}^{\ast}$.
\end{theorem}

\begin{remark}
Some versions of the semilocal splitting theorem and corresponding
counterexamples can also be found in~\cite{Fe1}.
\end{remark}

The natural question is to give some examples where condition
(4.28) does not hold and hence the Poisson structure is
nonlinearizable. First such examples of nonlinearizable 
Poisson structures at a symplectic leaf of nonzero dimension
were given in~\cite{DW}. These structures are constructed by
using a special class of Poisson structures called
\textit{Casimir-weighted products}.  We illustrate our results
with the following example.

\begin{example}
Let $B$ be an orientable $2$-surface with symplectic form
$$
\omega=ds\wedge d\tau.
$$
Consider the direct product $E=B\times\mathbb{R}^{3}$ as a
vector bundle over $B$, where $\pi$ is the canonical projection
on the first factor. Let 
$\mathbf{x}=(x_{1},x_{2},x_{3})$ be the Euclidean coordinates on
$\mathbb{R}^{3}$. Suppose we are given a vector-valued $1$-form 
$\boldsymbol{\varrho}$ on~$B$:
$$
\boldsymbol{\varrho}=\boldsymbol{\varrho}^{(1)}ds
+\boldsymbol{\varrho}^{(2)}d\tau,
$$
where $\boldsymbol{\varrho}^{(1)},\boldsymbol{\varrho}^{(2)}
\in\mathbb{R}^{3}$. Fix $c>0$. Then, one can associate to the
pair $(\boldsymbol{\varrho},c)$ the Poisson  
tensor $\Pi$ on $E$ which is defined by the following bracket relations:
\begin{align}
\{s,\tau\}&=\frac{1}{F(x)},
\tag{4.29}
\\
\{s,x_{\alpha}\}&=\frac{1}{F(x)}\epsilon_{\alpha\beta\gamma}x_{\beta}
\varrho_{\gamma}^{(2)},\qquad  
\{\tau,x_{\alpha}\}=-\frac{1}{F(x)}
\epsilon_{\alpha\beta\gamma} x_{\beta}\varrho_{\gamma}^{(1)},
\tag{4.30}
\\
\{x_{\alpha},x_{\beta}\}&=\epsilon_{\alpha\beta\gamma}x_{\gamma}+\frac{1}
{F(x)}\epsilon_{\alpha\sigma\gamma}\epsilon_{\beta\sigma^{\prime}
\gamma^{\prime}}[\varrho_{\sigma}^{(1)}\varrho_{\sigma^{\prime}}^{(2)}
-\varrho_{\sigma^{\prime}}^{(1)}\varrho_{\sigma}^{(2)}]x_{\gamma}
x_{\gamma^{\prime}}.
\tag{4.31}
\end{align}
Here $\epsilon_{ijk}$ is the completely antisymmetric
Levi-Civit\`a tensor,
$$
F(x)=1-\langle\mathbf{x},\mathbf{a}\rangle
+c\|\mathbf{x}\|^{2}
$$
and
$$
\mathbf{a}=\boldsymbol{\varrho}^{(1)}\times\boldsymbol{\varrho}^{(2)}.
$$
If
$$
\|\mathbf{a}\|^{2}<4c,
$$
then the brackets (4.29)--(4.30) are well defined on the entire
space $E$. It is clear that $(B,\omega)$ is the symplectic leaf
of this bracket. Using formulas (2.28)--(2.30), we compute the
corresponding geometric data: 
\begin{align}
\Pi_{V}&=\frac{1}{2}\epsilon_{\alpha\beta\gamma}x_{\gamma}
\frac{\partial}{\partial x^{\alpha}}\wedge
\frac{\partial}{\partial x^{\beta}},
\tag{4.32}
\\
\Gamma&=d\mathbf{x}-\mathbf{x}\times\boldsymbol{\varrho,}
\tag{4.33}
\\
\pi^{\ast}F&=F(x)ds\wedge d\tau.
\tag{4.34}
\end{align}
The curvature form of $\Gamma$ is
$$
Curv^{\Gamma}=\big\langle\mathbf{x}\times\mathbf{a},
\frac{\partial}{\partial\mathbf{x}}\big\rangle\otimes 
ds\wedge d\tau.
$$
The geometric data of the linearized Poisson structure
$\Pi^{(1)}$ are given by
\begin{align}
\Lambda&=\Pi_{V},
\tag{4.35}
\\
\Gamma^{(1)}&=\Gamma,
\tag{4.36}
\\
\pi^{\ast}F^{(1)}&=F^{(1)}(x)ds\wedge d\tau,
\tag{4.37}
\end{align}
where
\begin{equation}
F^{(1)}(x)=1-\langle\mathbf{x},\mathbf{a}\rangle.
\tag{4.38}
\end{equation}
One can recognize formula (4.32) as the Lie--Poisson structure
of the dual $\so^{\ast}(3)$ of the Lie algebra $\so(3)$. Thus,
the transverse Lie algebra of $B$ is $\so(3)$. It is clear that
$\Lambda$ is globally trivial. Remark that if $\mathbf{a}=0$,
then $\Pi$ is the Casimir-weighted product~\cite{DW}. Now,
comparing formulas (4.34) and (4.37), we compute the Casimir
$2$-cocycle $C$, 
\begin{equation}
\pi^{\ast}C=c\|\mathbf{x}\|^{2}ds\wedge d\tau.
\tag{4.39}
\end{equation}
The space of the Casimir functions of $\so^{\ast}(3)$ is
naturally identified with $C^{\infty}([0,\infty))$. Under this
identification, the coboundary operator is
$\partial_{0}=d_{B}\otimes1$, where $d_{B}$ is the exterior
differential on $B$. We conclude that 
if $B$ is compact, for example, $B$ is a sphere or a torus, 
then $[\omega]\neq0$ and hence $[C]\neq0$. 
By Theorem~4.14, in this case, the Poisson structure $\Pi$ is
nonlinearizable. On the other hand, if $\pi_{2}(B)=0$, for
example, $B=\mathbb{S}^{1} \times\mathbb{R}$ is the
$2$-cylinder, then $[C]=0$. In this case, $\Pi$ is linearizable
and equivalent to the direct product Poisson structure on 
$B\times\so^{\ast}(3)$.
\end{example}


\begin{thebibliography}{CrFe}
\bibitem[Co$_{1}$]{Co1}
J.~Conn, 
\textit{Normal forms for analytic Poisson structures}, 
Ann. Math., \textbf{119} (1984), 576--601.

\bibitem[Co$_{2}$]{Co2}
J.~Conn, 
\textit{Normal forms for smooth Poisson structures}, 
Ann. Math., \textbf{121} (1985), 565--593.

\bibitem[DW]{DW}
B.~L.~Davis and A.~Wade, 
\textit{Nonlinearizability of Certain Poisson Structures 
near a Symplectic Leaf}, 
Preprint, math.SG/0406611, 2004.

\bibitem[Vo$_{1}$]{Vo1}
Yu.~Vorobjev, 
\textit{Coupling tensors and Poisson geometry near a single
symplectic leaf}. 
In: \textit{Lie Algebroids}, 
Banach Center Publ., Vol.~54, Waszawa, 2001, pp.~249--274.

\bibitem[Vo$_{2}$]{Vo2}
Yu.~Vorobjev, 
\textit{On linearized Poisson brackets},
Math. Notes, \textbf{70} (2001), no.~4, 486--493.

\bibitem[Vo$_{3}$]{Vo3}
Yu.~Vorobjev, 
\textit{On Poisson realizations of transitive Lie algebroids}, 
J. Nonlinear Math. Phys., \textbf{11} (2004), 43--48.

\bibitem[Vo$_{4}$]{Vo4}
Yu.~Vorobjev, 
\textit{Poisson structures and linear Euler systems over
symplectic manifolds}, 
Amer. Math. Soc. Transl. (2), Providence, RI, 
in preparation.

\bibitem[IKV]{IKV}
V.~M.~Itskov, M. Karasev, and Yu.~M.~Vorobjev,
\textit{Infinitesimal Poisson cohomology}, 
Amer. Math. Soc. Transl. (2), Vol.~187, 
Providence, RI, 1998, 327--360.

\bibitem[GLS]{GLS}
V.~Guillemin, E.~Lerman, and S.~Sternberg,
\textit{Symplectic Fibrations and Multiplicity Diagrams}, 
Cambridge Univ. Press., Cambridge, 1996.

\bibitem[We$_{1}$]{We1}
A.~Weinstein, 
\textit{The local structure of Poisson manifolds}, 
J. Diff. Geom., \textbf{18} (1983), 523--557.

\bibitem[We$_{2}$]{We2}
A.~Weinstein, 
\textit{Symplectic groupoids and Poisson manifolds}, 
Bull. Amer. Math. Soc., \textbf{16} (1987), 101--104.

\bibitem[We$_{3}$]{We3}
A.~Weinstein, 
\textit{Linearization problem for Lie algebroids and Lie groupoids}, 
Lett. Math. Phys., \textbf{52} (1987), 55--73.

\bibitem[Br]{Br}
O.~Brahic, 
\textit{Normal Forms of Poisson Structures Near a Symplectic Leaf}, 
Preprint, math.SG/0403136, 2004.

\bibitem[Fe$_{1}$]{Fe1}
R.~L.~Fernandes, 
\textit{Connections in Poisson geometry I: 
holonomy and invariants}, 
J. Differential Geom., \textbf{54} (2000), no.~2, 303--365.

\bibitem[Fe$_{2}$]{Fe2}
R.~L.~Fernandes, 
\textit{Lie algebroids, holonomy, and characteristic classes},
Adv. in Math., \textbf{170} (2002), 119--179.  

\bibitem[FeM]{FeM}
R.~L.~Fernandes and P.~Monnier, 
\textit{Linearization of Poisson Brackets}, 
Preprint, math.SG/ 0401273, 2004.

\bibitem[Va$_{1}$]{Va1}
I.~Vaisman, 
\textit{Lectures on the Geometry of Poisson Manifolds}, 
Progress in Math., Vol.~118, Birkh\"auser, Boston, 1994.

\bibitem[Va$_{2}$]{Va2}
I.~Vaisman, 
\textit{Coupling Poisson and Jacobi structures on foliated
manifolds},  
J. of Geometric Methods in Modern Physics, 
\textbf{1} (2004), no.~5, 607--637.

\bibitem[Va$_{3}$]{Va3}
I.~Vaisman, 
\textit{Foliation-Coupling Dirac Structures}, 
Preprint, math.SG/0412318, 2004.

\bibitem[Ar]{Ar}
V.~I.~Arnold, 
\textit{Mathematical Methods of Classical Mechanics}, 
Graduate Text in Math., Vol.~60, Springer-Verlag, 
New York, 1978.

\bibitem[DZ]{DZ}
J.-P.~Dufour and N.~T.~Zung, 
\textit{Nondegeneracy of Lie algebra $\mathfrak{aff}\,(n)$}, 
C. R. Acad. Sci. Paris S\'{e}r. I Math., 
\textbf{335} (2002), 1043--1046.

\bibitem[Du$_{1}$]{Du1}
J.-P.~Dufour, 
\textit{Lin\'{e}arisation de Certaines Structures de Poisson}, 
J. Diff. Geom., \textbf{32} (1990), no.~2, 415--428.

\bibitem[Du$_{2}$]{Du2}
J.-P.~Dufour, \textit{Normal forms of Lie algebroids},
Banach Center Publ., \textbf{54} (2001), 35--41.

\bibitem[Zu$_{1}$]{Zu1}
N.~T.~Zung, 
\textit{Levi decomposition of analytic Poisson structures and
Lie algebroids}, 
Topology, \textbf{42} (2003), no.~6, 1403--1420. 

\bibitem[Zu$_{2}$]{Zu2}
N.~T.~Zung, 
\textit{A geometric Proof of Conn's Linearization Theorem for
Analytic Poisson Structures}, 
Preprint, math.SG/0207263, 2002. 

\bibitem[St]{St}
S.~Sternberg, 
\textit{Minimal coupling and the symplectic mechanics of a
classical particle in the presence of a Young--Mills field}, 
Proc. Nat. Acad. Sci. U.S.A., \textbf{74} (1977), 5253--5254.

\bibitem[KM]{KM}
M.~V.~Karasev and V.~P.~Maslov, 
\textit{Nonlinear Poisson Brackets. Geometry and Quantization}, 
Transl. of Math. Monographs, Vol.~119,
Amer. Math. Sos., Providence, RI, 1993.

\bibitem[GHV]{GHV}
W.~Greub, S.~Halperin, and R.~Vanstone, 
\textit{Connections, Curvature, and Cohomology}, 
Vol.~II, Academic Press, New York--London, 1973.

\bibitem[Ku]{Ku}
J.~Kubarski, 
\textit{The Chern-Weil homomorphism of regular Lie algebroids}, 
Publ. Dep. Math. Nouvelle Ser., 
Univ. Claude-Bernard Lyon 1, 1991, 1--69.

\bibitem[Mz]{Mz}
K.~C.~H.~Mackenzie, 
\textit{Lie Groupoids and Lie Algebroids in Differential Geometry}, 
LMS Lecture Note Ser., Vol.~124, 
Cambridge Univ. Press, Cambridge, 1987.

\bibitem[LMr]{LMr}
P.~Libermann and C.-M.~Marle, 
\textit{Symplectic Geometry and Analytical Mechanics}, 
Reidel, Dordrecht, 1987.

\bibitem[GiGo]{GiGo}
V.~L.~Ginzburg and A.~Golubev, 
\textit{Holonomy on Poisson manifolds and the modular class}, 
Israel. J. Math., \textbf{122} (2001), 221--242. 

\bibitem[CrFe]{CrFe} 
M.~Crainic and R.~L.~Fernandes, 
\textit{Rigidity and Flexibility in Poisson Geometry}, 
Preprint math.DG/0503145 v1, 2005.

\end{thebibliography}
\end{document}